\documentclass[twoside,11pt]{article}

\usepackage{graphicx, subfigure}
\usepackage[top=1in,bottom=1in,left=1in,right=1in]{geometry}
\usepackage[sort&compress]{natbib} 
\usepackage{amsfonts, amsmath, amssymb, amsthm, constants, bbm}
\usepackage{MnSymbol}
\usepackage{mathtools}
\usepackage[mathscr]{eucal}

\usepackage{color}
\definecolor{darkred}{RGB}{100,0,0}
\definecolor{darkgreen}{RGB}{0,100,0}
\definecolor{darkblue}{RGB}{0,0,150}

\usepackage{hyperref}
\hypersetup{colorlinks=true, linkcolor=darkred, citecolor=darkgreen, urlcolor=darkblue}
\usepackage{url}

\newtheorem{thm}{Theorem}

\newtheorem{lem}{Lemma}
\newtheorem{cor}{Corollary}

\theoremstyle{remark}
\newtheorem{Def}{Definition}
\newtheorem{rem}{Remark}

\def\beq{\begin{equation}} 
\def\eeq{\end{equation}}
\def\beqn{\begin{eqnarray*}}
\def\eeqn{\end{eqnarray*}}
\def\Bitem{\begin{itemize}\setlength{\itemsep}{.2in}}
\def\bitem{\begin{itemize}\setlength{\itemsep}{.05in}}
\def\eitem{\end{itemize}}
\def\Benum{\begin{enumerate}\setlength{\itemsep}{.2in}}
\def\benum{\begin{enumerate}\setlength{\itemsep}{.05in}}
\def\eenum{\end{enumerate}}
\def\bmult{\begin{multline*}}
\def\emult{\end{multline*}}
\def\bcenter{\begin{center}}
\def\ecenter{\end{center}}
\def\bframe{\begin{frame}}
\def\eframe{\end{frame}}

\newcommand{\thmref}[1]{Theorem~\ref{thm:#1}}

\newcommand{\corref}[1]{Corollary~\ref{cor:#1}}
\newcommand{\lemref}[1]{Lemma~\ref{lem:#1}}
\newcommand{\secref}[1]{Section~\ref{sec:#1}}
\newcommand{\figref}[1]{Figure~\ref{fig:#1}}

\newcommand{\remref}[1]{Remark~\ref{rem:#1}}
\newcommand{\defref}[1]{Definition~\ref{def:#1}}

\DeclareMathOperator*{\argmin}{arg\, min}

\def\cF{\mathcal{F}}

\def\cN{\mathcal{N}}

\def\bX{\mathbf{X}}

\newcommand\bPsi{{\boldsymbol\Psi}}

\def\bbH{\mathbb{H}}
\def\bbI{\mathbb{I}}

\def\bbP{\mathbb{P}}

\def\bbR{\mathbb{R}}

\newcommand{\E}{\operatorname{\mathbb{E}}}
\renewcommand{\P}{\operatorname{\mathbb{P}}}

\newcommand\indep{\protect\mathpalette{\protect\independenT}{\perp}}
\def\independenT#1#2{\mathrel{\rlap{$#1#2$}\mkern2mu{#1#2}}}

\def\Bin{\text{Bin}}

\def\eps{\varepsilon}

\def\1{\mathbbm{1}}
\newcommand{\IND}[1]{\bbI\{ #1 \}}

\def\A{\mathscr{A}}
\def\F{\mathscr{F}}

\def\R{\mathscr{R}}

\def\H{\bbH}

\def\fdr{\textsc{fdr}}
\def\fdp{\textsc{fdp}}
\def\fnr{\textsc{fnr}}
\def\fnp{\textsc{fnp}}
\def\fdphat{\widehat{\textsc{fdp}}}

\definecolor{purple}{rgb}{0.4,.1,.9}

\newcommand\blfootnote[1]{%
  \begingroup
  \renewcommand\thefootnote{}\footnote{#1}%
  \addtocounter{footnote}{-1}%
  \endgroup
}

\begin{document}
\thispagestyle{empty}

\title{Distribution-free Multiple Testing}
\author{
Ery Arias-Castro 
\and Shiyun Chen
}
\date{}
\maketitle

\begin{abstract}
We study a stylized multiple testing problem where the test statistics are independent and assumed to have the same distribution under their respective null hypotheses.  We first show that, in the normal means model where the test statistics are normal Z-scores, the well-known method of \citep{benjamini1995controlling} is optimal in some asymptotic sense.  We then show that this is also the case of a recent distribution-free method proposed by \cite{barber2015controlling}.  The method is distribution-free in the sense that it is agnostic to the null distribution --- it only requires that the null distribution be symmetric.  We extend these optimality results to other location models with a base distribution having fast-decaying tails. 
\end{abstract}

\blfootnote{Both authors are with the Department of Mathematics, University of California, San Diego, USA.  Contact information is available \href{http://math.ucsd.edu/\~eariasca}{here} and \href{http://www.math.ucsd.edu/people/graduate-students/}{here}.  This work was partially supported by grants from the US National Science Foundation (DMS 1223137) and the US Office of Naval Research (N00014-13-1-0257).}  


\section{Introduction} \label{sec:intro}
Multiple testing arises in a wide array of applied settings, ranging from anomaly detection in sensor arrays to the selection of genes that are differentially expressed \citep{dudoit2007multiple, dickhaus2014simultaneous}.  This is particularly true in so-called discovery science, where the scientist proceeds by formulating hypotheses, testing each one of them on data, and following up on the most promising ones.  Each step along the way is fraught with pitfalls, yet even if the experiment was correctly designed and carried out, the scientist still needs to contend with the multitude of tests that were performed.  

Multiple testing is now a well-established area in statistics.  In a substantial proportion of the corresponding literature it is assumed that P-values are available.  This, implicitly, assumes that the null distribution of each test statistic is known (perfectly).  For example, the Benjamini-Hochberg (BH) procedure was proposed in this context \citep{benjamini1995controlling}.  
See \citep{roquain2011type} for a fairly recent and comprehensive review of the literature, as it pertains to mathematical results in the area.

Our contribution is two-fold.  First, we prove that the BH method is asymptotically optimal to first order in the normal (location) model, which corresponds to an idealized setting where the tests being performed are Z-tests and the effect, when present, affects the mean.  In fact, we show that this is the case in the much wider context of asymptotically generalized Gaussian models --- see \defref{AGG}.  
Second, we propose to use the recent distribution-free method of \cite{barber2015controlling} that only relies on the assumption that the test statistics have a common null distribution that is symmetric about 0 and show that, in the same normal model, it achieves the same asymptotic performance to first order.  This method, proposed in the context of post-model selection inference, is also intimately related to our own work \citep{arias2013distribution} on distribution-free testing of the global null hypothesis.

\subsection{The risk of a multiple testing procedure}
Consider a setting where we want to test $n$ null hypotheses, denoted $\H_1, \dots, \H_n$.  The test that we use for $\H_i$ rejects for large values of a statistic $X_i$.  Throughout, we assume that $X_1, \dots, X_n$ are independent.  Denote the vector of test statistics by $\bX = (X_1, \dots, X_n)$. 
Let $\Psi_i$ denote the survival function\footnote{In this paper, the survival function of a random variable $Y$ is defined as $y \mapsto \P(Y \ge y)$.} of $X_i$ and $\bPsi = (\Psi_1, \dots, \Psi_n)$.

\begin{rem}
In a large portion of the literature, it is assumed that P-values can be computed (or at least approximated).  The simplest such case is when $\H_i$ is a singleton, $\H_i = \{\Psi_i^{\rm null}\}$, and the null distributions $\Psi_1^{\rm null}, \dots, \Psi_n^{\rm null}$ are known.  In that case, the $i$-th P-value is defined as $P_i = \Psi_i^{\rm null}(X_i)$, which is the probability of exceeding the observed value of the statistic under its null distribution.  In this context, working with the statistics $X_1, \dots, X_n$ is equivalent to working with the P-values $P_1,\dots,P_n$.
\end{rem}

Let $\F \subset [n] := \{1, \dots, n\}$ index the false null hypotheses, meaning
\beq
\F = \{i \in [n] : \Psi_i \notin \H_i\}.
\eeq
A multiple testing procedure $\R$ takes the test statistics $\bX$ and return a subset of $\R(\bX) \subset [n]$ representing the null hypotheses that the procedure rejects.  Given such a procedure $\R$, the false discovery rate is defined as the expected value of the false discovery proportion \citep{benjamini1995controlling}
\beq
\fdr_\bPsi(\R) = \E_\bPsi(\fdp(\R(\bX))), \quad \fdp(\R) := \frac{|\R \setminus \F|}{|\R|},
\eeq
where we denoted the cardinality of a set $\A \subset [n]$ by $|\A|$ and we adopt the convention that $0/0 = 0$.
While the FDR of a multiple testing procedure is analogous to the level or size of a test procedure, the false non-discovery rate plays the role of power and is defined as the expected value of the false non-discovery proportion
\beq
\fnr_\bPsi(\R) = \E_\bPsi(\fnp(\R(\bX))), \quad \fnp(\R) := \frac{|\F \setminus \R|}{|\F|}.
\eeq
(Note that our definition is different from that of \cite{genovese2002operating}.  )
Accordingly, we define the risk of a multiple testing procedure $\R$ as
\beq\label{risk}
{\rm risk}_\bPsi(\R) = \fdr_\bPsi(\R) + \fnr_\bPsi(\R).
\eeq

\begin{rem} \label{rem:fnp-fnr}
For any multiple testing procedure, $\fnp \to 0$ in probability if and only if $\fnr \to 0$.  Indeed, one direction is justified by Markov's inequality, and the other direction is justified by dominated convergence (using the fact that $\fnp \le 1$ always).
\end{rem}

\subsection{Natural procedures}
Under our circumstances, the following properties are {\em natural}, and are for example satisfied by the BH procedure, and in fact, any other procedure we can think of, at least in the present context.

\bitem
\item {\em Permutation invariance.} 
We say that a procedure $\R$ is permutation invariant if, as a function of $n$ variables, it is invariant with respect to the order of these variables.  \item {\em Threshold.}
We say that a procedure $\R$ is a threshold procedure if it is of the form
\beq\label{threshold}
\R(X_1, \dots, X_n) = \{i : X_i \ge \tau(X_1, \dots, X_n)\},
\eeq
for some threshold function $\tau$.
\item {\em Monotonicity.}
We say that a procedure $\R$ is monotonic if, for any $i \in \R(X_1, \dots, X_n)$, increasing $X_i$ leaves $\R$ unchanged. 
\eitem

Because they are so natural, but also for ease of exposition, we focus on natural procedures throughout the paper.


\subsection{The normal model and the optimality of the BH method}
This model corresponds to the setting above with $X_i \sim \Psi_i = \cN(\mu_i, 1)$ and $\H_i : \mu_i = 0$, so that $\H_i$ is a singleton equal to $\Psi_i^{\rm null} = \cN(0,1)$.  In this context it is compelling to ask how the $\mu_i$'s need to be in order for the risk of the BH procedure to tend to zero.  To the best of our knowledge, this question has not been directly answered in the literature.

Our inspiration for considering the normal (location) model comes from the seminal work of Ingster \citep{ingster1997some,IngsterBook} and \cite{donoho2004higher} on testing the global null $\bigcap_i \H_i$.  
In \citep{ingster1997some} we  find the following first-order asymptotic result.  Assume a prior under which $m \le n$ randomly picked $\mu_i$'s are set to $\sqrt{2 r \log n}$ and the others are set to 0.  An interesting parameterization happens to be $m/n \sim n^{-\beta}$ with $\beta > 0$ fixed.  Focusing on the so-called sparse regime, where $\beta > 1/2$, one finds that the detection boundary is at $r = \rho(\beta)$, where 
\beq \label{lower}
\rho(\beta) = \begin{cases} 
\beta - 1/2, & 1/2 < \beta \leq 3/4, \\ 
(1 - \sqrt{1 - \beta})^2, & 3/4 < \beta < 1.
\end{cases}
\eeq
This means that, taking $r$ to be fixed, when $r < \rho(\beta)$ all tests have risk\footnote{As usual, the risk of a test is defined as the sum of its type I and type II error probabilities.} at least 1 in the large sample limit (which is as bad as random guessing), while when $r > \rho(\beta)$ the  likelihood ratio test has risk 0 in the large sample limit.
\cite{donoho2004higher} propose an adaptive test procedure based on Tukey's higher criticism (HC) that achieves this optimal detection boundary.   

Returning to the question of identifying the false null hypotheses, which is our concern here, we know that $r > 1$ allows for the identification of the false nulls with a control of the family-wise error rate (FWER) at any fixed level.  In fact, $r = 1$ is the precise boundary for this to be possible --- we leave this as an exercise to the reader.  This is true for all $\beta \in (0,1)$.
We are more interested in controlling the risk \eqref{risk}.  The following is a special case of a more general lower bound appearing later in the paper.

\begin{cor} \label{cor:normal-lb}
In the normal model, assume that $\beta \in (0,1)$ and $r \ge 0$ are both fixed.
If $r < \beta$, then the risk of any natural multiple testing procedure has limit inferior at least 1 as $n \to \infty$.  
\end{cor}

\begin{rem}
Such a lower bound was in fact obtained earlier by \cite{ji2012ups} in a more general normal model allowing for dependencies.  
\end{rem}

In our context, we know that \corref{normal-lb} is tight because the BH method 
achieves the stated multiple testing boundary.  The following is also a special case of a more general result appearing later in the paper.

\begin{cor} \label{cor:normal-BH}
In the setting of \corref{normal-lb}, if instead $r > \beta$, then the risk of the BH procedure (properly calibrated) tends to 0 as $n \to \infty$.  
\end{cor}

\begin{rem}
A similar result was derived in \citep{jin2014rare} for a form of oracle threshold method with threshold chosen based on knowledge of the number of false nulls.  By comparison, the BH method --- which does not require this knowledge --- is shown here to be adaptive to the number of false nulls.    
\end{rem}

Together, \corref{normal-lb} and \corref{normal-BH} establish that the BH procedure is asymptotically optimal to first order in the normal model.  We will see that this remains true for a much wider class of models. 


\subsection{Multiple testing under symmetry}
The P-values are based on the assumed knowledge of the null distribution of each test statistic.  In many practical settings, this is not strictly the case, resulting in P-values that are only approximately uniformly distributed under their respective null hypothesis.  This can jeopardize the control of the FDR.  In the same way that it may be appealing in some situations to use a distribution-free test such as the signed-rank test instead of the t-test, it may be desirable to use a distribution-free procedure for multiple testing.  

Our working assumption is the following
\beq\label{sym}
\textit{$X_1, \dots, X_n$ are independent with common null distribution that is symmetric about 0.}
\eeq
The assumption of symmetry is standard in the literature on nonparametric tests \citep{MR758442}.  However, it is not standard in the context of multiple testing.  Although quite natural, we only know of two instances where this assumption is made:
\bitem
\item In our own work on testing the global null \citep{arias2013distribution}.  
\item In a method for post-Lasso model selection inference \citep{barber2015controlling}.
\eitem 
Whether following one lead or the other, the resulting procedure is essentially the same, which we call the Barber-Cand\`es (BC) procedure, properly defined later on.  
This procedure is shown in \citep{barber2015controlling} to control the FDR at any desired level.  
Here we show that, under fairly general conditions, it achieves the multiple testing boundary.  In particular, it does as well as the BH procedure with knowledge of the null distributions.  
The following is a special case of a more general results appearing later on.

\begin{cor} \label{cor:normal-BC}
The conclusions of \corref{normal-BH} apply to the BC procedure.
\end{cor}

We mention that other distribution-free procedures have been suggested in the literature.  Except for those cited above, all the other ones we know of are based on resampling \citep{westfall1993resampling,
ge2003resampling,
yekutieli1999resampling,
romano2007control}.
These methods are not applicable in the setting assumed here.  They are typically applied to situations, as in microarray analysis, where each test statistic is based on comparing two (or more) samples.

Another approach suggested in the literature is that of estimating the null distribution of the test statistics (assumed to be the same).  This is advocated in \citep{efron2004large,
pollard2004choice}.  The asymptotic risk properties of such methods remain unknown, despite the fact that there is some theory on the topic of estimating the null distribution \citep{jin2007estimating,
cai2010optimal}.  

\subsection{Content}
In \secref{lower} we derive a lower bound on the boundary for multiple testing in a location model where the base distribution is asymptotically generalized Gaussian.  This comprises the normal model.
In \secref{BH} we analyze the performance of the BH procedure based on the full knowledge of the null distribution, while in \secref{BC} we analyze the performance of the BC procedure.
We present the result of some numerical experiments in \secref{numerics}.
The proofs are gathered in \secref{proofs}.

\section{The AGG model}
\label{sec:lower}

We start by defining an oracle procedure, which is all we can hope for when using a natural procedure.  We then define a family of location models where the base distribution is asymptotically polynomial in log-scale --- which in particular encompasses the normal model --- and in the context of such a model, we establish a lower bound on the risk of the oracle procedure.  The result is an oracle risk bound.

\subsection{The oracle procedure}
In the context of natural procedures, the oracle procedure is defined as the threshold procedure using the threshold that minimizes the risk of a particular realization, namely,
\beq\label{oracle}
\tau_o \in \argmin_{t \in \bbR} ~ \fdp(\R_t) + \fnp(\R_t), \qquad \R_t := \{i : X_i \ge t\}.
\eeq
In words, with full knowledge of the set of false null distributions $\F$, the procedure chooses a threshold that partitions the test statistics in a way that minimizes the sum of the false discovery and non-discovery proportions.  The expected risk of this procedure is what we call below the oracle risk.

\begin{rem}
Of course, if one knew $\F$, one would simply reject $\H_i$ for all $i \in \F$ and, in the end, there would not any multiple testing problem to deal with!  This oracle procedure is, however, constrained to be of threshold type.  We will only use this fictitious procedure as a benchmark among threshold-type procedures. 
\end{rem}

\begin{rem}
We note that our oracle is stronger (provides more information) than the oracle considered by \cite{meinshausen2011asymptotic} in the context of FWER control.
\end{rem}

\subsection{Asymptotically generalized Gaussian model} \label{sec:model}
In a location model, we assume that we know the null survival function $\Psi$, assumed to be continuous for simplicity, and the test statistics are independent with respective distribution $X_i \sim \Psi_i = \Psi(\cdot - \mu_i)$, where $\mu_i = 0$ under the null $\H_i$ and $\mu_i > 0$ otherwise.
Both minimax and Bayesian considerations lead to considering a prior on the $\mu_i$'s where $m \le n$ randomly picked $\mu_i$'s are set equal to some $\mu > 0$ and the others are set to 0.  The prior is therefore defined based on $m$ and $\mu$, which together control the signal strength.

Beyond the normal model, we consider other location models where the base distribution has a polynomial right tail in log scale.  

\begin{Def} \label{def:AGG}
A survival function $\Psi$ is asymptotically generalized Gaussian (AGG) on the right with exponent $\gamma > 0$ if $\lim_{x \to \infty} x^{-\gamma} \log \Psi(x) = -1/\gamma$.
\end{Def}

The AGG class of distributions is nonparametric and quite general.  
It includes the parametric class of generalized Gaussian (GG) distributions with densities $\{\psi_\gamma, \gamma > 0\}$ given by $\log \psi_\gamma(x) \propto -|x|^\gamma/\gamma$, which comprises the normal distribution ($\gamma = 2$) and the double exponential distribution ($\gamma=1$).  We assume that $\gamma \ge 1$ so that the null distribution has indeed a sub-exponential right tail.

\begin{rem}
We note that the scale (e.g., standard deviation) is fixed, but this is really without loss of generality as both the BH and BC methods are scale invariant.  For the BH method, this is because the P-values are scale invariant.  However, this is so because we provide the BH method with the null distribution, including the scale.  The BC method, by contrast, can operate without knowledge of the scale.
\end{rem}

\citep{donoho2004higher} considered the problem of testing the global null in a GG location model and derived the detection boundary.  
We use the same prior, where $m$ nulls chosen uniformly at random are designated to be false and all positive $\mu_i$'s are set equal to $\mu$, with
\beq \label{eps}
m = \lfloor n^{1-\beta} \rfloor, \quad \text{ with } 0 < \beta < 1 \quad \text{(fixed)},
\eeq  
and
\beq \label{mu}
\mu = \mu_\gamma(r)= \left( \gamma r \log n  \right)^{1/\gamma}, \quad \text{ with } r > 0 \quad \text{(fixed)}.
\eeq

\begin{thm} \label{thm:lower}
Consider a location model where the base distribution is AGG with exponent $\gamma \ge 1$, with prior described above, and with the parameterization \eqref{eps}-\eqref{mu}.
If $r < \beta$, then the oracle risk has limit inferior at least 1 as $n \to \infty$.  
\end{thm}

\begin{rem}
We believe the result is valid for all procedures, regardless of whether they are natural or not.  As usual, proving such a result would require the application of tools from information theory or decision theory  \citep{MR2724359}.  Fano's inequality is appealing, but to obtain a result as strong as \thmref{lower}, one would need to develop a stronger form of Fano's inequality.  
\end{rem}

\section{The performance of the BH method}
\label{sec:BH}

Recall that $P_i = \Psi(X_i)$ is the P-value corresponding to the test statistic $X_i$, where $\Psi$ denotes the survival function under the null.  
We order the $X_i$'s in {\em decreasing} order, to obtain the following order statistics $X_{(1)} \ge \cdots \ge X_{(n)}$.  
Given a desired FDR control at $q$, the BH procedure of \citep{benjamini1995controlling} is defined is the threshold procedure \eqref{threshold}, with threshold
\beq \label{def_BH}
\tau_{\rm BH} = X_{(\iota_{\rm BH})}, \quad 
\iota_{\rm BH} := \max\big\{i : P_{(i)} \le i q/n\big\},
\eeq
where $P_{(1)} \le \cdots \le P_{(n)}$ are the ordered P-values.
This procedure is shown in \citep{benjamini1995controlling} to control the FDR at $q$, for example, when the tests are independent --- which we assume throughout.

Typically, $q$ is set to a small number, like $q = 0.10$.  In this paper we allow $q$ to tend to 0 as $n \to \infty$, but slowly.  Specifically, we always assume that
\beq\label{q}
\text{$q>0$ such that $n^a q \to \infty$ for any $a > 0$ fixed.}
\eeq

The following result establishes the BH procedure as optimal in the AGG model, in the sense that it achieves the detection boundary ($r = \beta$) stated in \thmref{lower}.

\begin{thm} \label{thm:BH}
In the setting of \thmref{lower}, if instead $r > \beta$, then the risk of the BH procedure with $q$ satisfying \eqref{q} has $\fnp$ tending to 0 in probability as $n \to \infty$.  In particular, if $q \to 0$, then its risk tends to 0. 
\end{thm}

We note that the second part of the theorem follows from \remref{fnp-fnr} and that \corref{normal-BH} follows immediately from this result.

\section{The performance of the BC method}
\label{sec:BC}

Under the assumption of symmetry, given the desired FDR control level $q$, the Barber-Cand\`es (BC) procedure defines the data-dependent threshold $\tau_{\rm BC}$ as:

\beq  \label{BC}
\tau_{\rm BC} = \inf\big\{t \in |\bX|: \fdphat(t) \leq q\big\},
\eeq
where, as usual, the infimum is infinite if the set is empty, $|\bX| := \{|X_i| : i = 1, \dots, n\}$ is the set of sample absolute values, and
\beq
\fdphat(t) := \frac{1+\#\{i : X_i \leq -t\}}{1 \vee \#\{i: X_i \geq t\}},
\eeq
is a measure of how asymmetric the set of observations $\{X_i: |X_i| \geq t\}$ is.  

The notation is borrowed from \citep{barber2015controlling} and is justified by the fact that this quantity aims at estimating $\fdp(\R_t)$, where $\R_t = \{i : X_i \ge t\}$ as in \eqref{oracle}.
The BC procedure is shown in \citep{barber2015controlling} to control the FDR at level $q$.

The following result shows that, although agnostic to the null distribution, the BC procedure achieves the detection boundary in a AGG model as long as the underlying distribution is symmetric.

\begin{thm} \label{thm:BC}
In the setting of \thmref{lower}, and assuming that the null distribution is symmetric about 0, if instead $r > \beta$, then the risk of the BC procedure with $q$ satisfying \eqref{q} has $\fnp$ tending to 0 in probability as $n \to \infty$.  In particular, if $q \to 0$, then its risk tends to 0. 
\end{thm}

\subsection{Another variant}
We mentioned another inspiration for considering the BC procedure as a potential candidate for achieving the detection boundary in such a nonparametric setting --- which is now confirmed by \thmref{BC}.  
It relates to our own work on testing the global null in a similar setting \citep{arias2013distribution}.
Following closely the reasoning there leads us to consider the following procedure.  Let $\xi_{(1)}, \dots, \xi_{(n)}$ be the respective signs of the observations arranged in decreasing order of absolute value and let $S_k = \sum_{i \le k} \xi_{(i)}$ denote their partial sum up to $k$.  Given a desired $q \in (0,1)$, we define the threshold index 
\beq
\iota = \max\big\{i : S_i \ge i \tfrac{1-q}{1+q}\big\}.
\eeq
The procedure is then defined as the threshold procedure with threshold $\tau = |X_{(\iota)}|$, where $|X|_{(i)}$ denotes the $i$-th largest observation in absolute value.
To make a strong parallel with \citep{arias2013distribution}, let us call this procedure the cumulative sum (CUSUM) sign procedure.  

The two procedures (BC and CUSUM sign) are equally principled, and in fact, they are very closely related, as we show next.
Assume that the observations come from a continuous distribution, so that no observation is equal to 0 with probability 1, or that the observations equal to 0 have been removed.  In that case,
\beq
S_k = \sum_{i=1}^{k} \xi_{(i)} = \sum_{i=1}^{k} \IND{\xi_{(i)} > 0} - \sum_{i=1}^{k} \IND{\xi_{(i)} < 0} = 2\sum_{i=1}^{k} \IND{\xi_{(i)} > 0}-k.
\eeq
We can express $\fdphat(|X|_{(k)})$ as a function of $S_k$:
\beq
\fdphat(|X|_{(k)}) = \frac{1+\#\{j: X_j \leq -|X|_{(k)}\}}{1 \vee \#\{j: X_j \geq |X|_{(k)}\}} 
= \frac{1+\sum_{i=1}^{k} \IND{\xi_{(i)} < 0}}{1 \vee \sum_{i=1}^{k} \IND{\xi_{(i)} > 0}}
= \frac{1 + (k - S_k)/2}{1 \vee (k + S_k)/2}.
\eeq
Then, under mild assumptions, 
\beq
\fdphat(|X|_{(k)}) \approx \frac{k - S_k}{k + S_k},
\eeq
and when this is the case the event $\fdphat(|X|_{(k)}) \le q$ is approximately equivalent to the event $S_k \ge (1-2q/(1+q)) k$.

Thus we suspect that the CUSUM sign procedure performs comparably to the BC procedure.  We did observed this in all the numerical experiments we performed (none reported here).

\section{Numerical experiments}
\label{sec:numerics}
In this section, we perform simple simulations to compare the BH and BC procedures on finite data, with the goal of illustrating the theory we established.  We consider the normal model and the double-exponential model.  
We reemphasize that the BH procedure plays a role of oracle here since it requires knowledge of null distribution to compute the P-values.  In contrast, the BC method does not require knowledge of the null distribution.

\subsection{Fixed sample size} \label{sec:fixed}
In this first set of experiments, the sample size is chosen large at $n = 10^5$. 
The FDR control level is set at $q = 0.05$. We draw $m$ observations from the alternative distribution $\Psi(\cdot - \mu)$, and the other $n-m$ from the null distribution $\Psi$.  All the models are parameterized as described in \secref{model}, in particular, \eqref{eps} and \eqref{mu}. 
We choose a few values for the parameter $\beta$ so as to exhibit different sparsity levels, while the parameter $r$ takes values in a grid of range $(0,1)$.  Each situation is repeated 100 times for each test and we report the average FDP and FNP of both tests.

\subsubsection{Normal model}
\begin{figure}[h!]\centering
	\includegraphics[width = 5cm, height= 5cm]{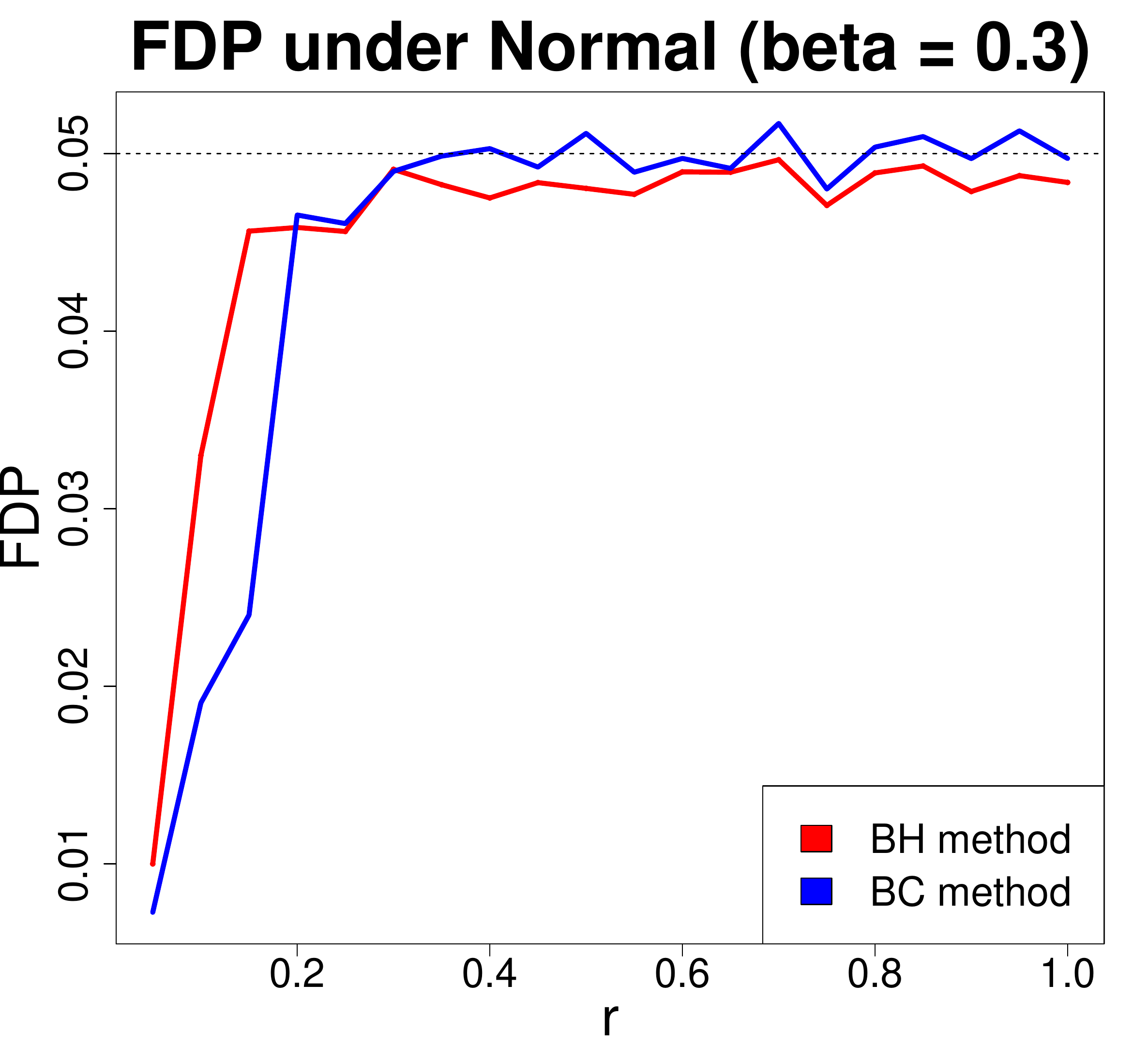}
	\includegraphics[width = 5cm, height= 5cm]{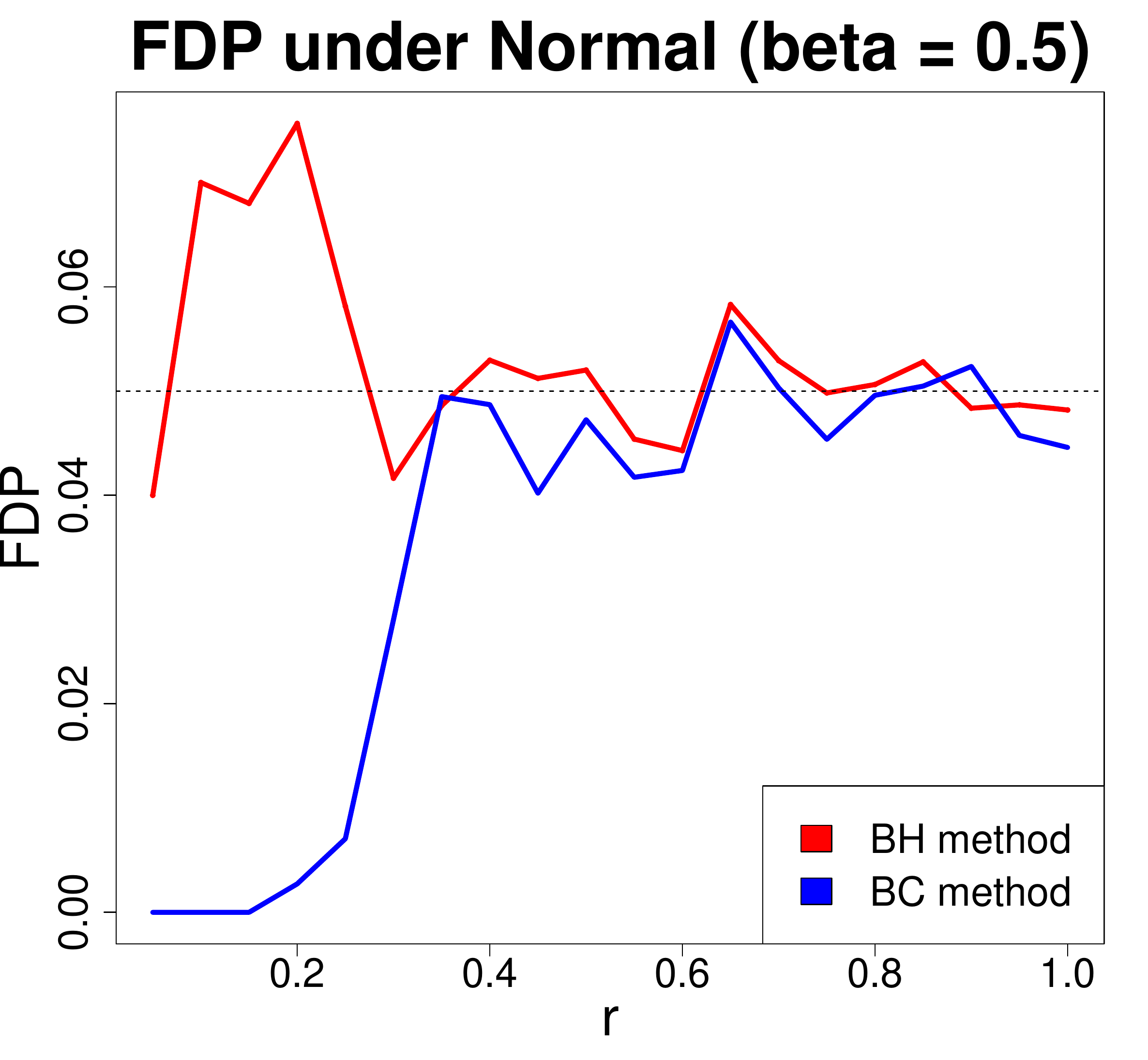}
	\includegraphics[width = 5cm, height= 5cm]{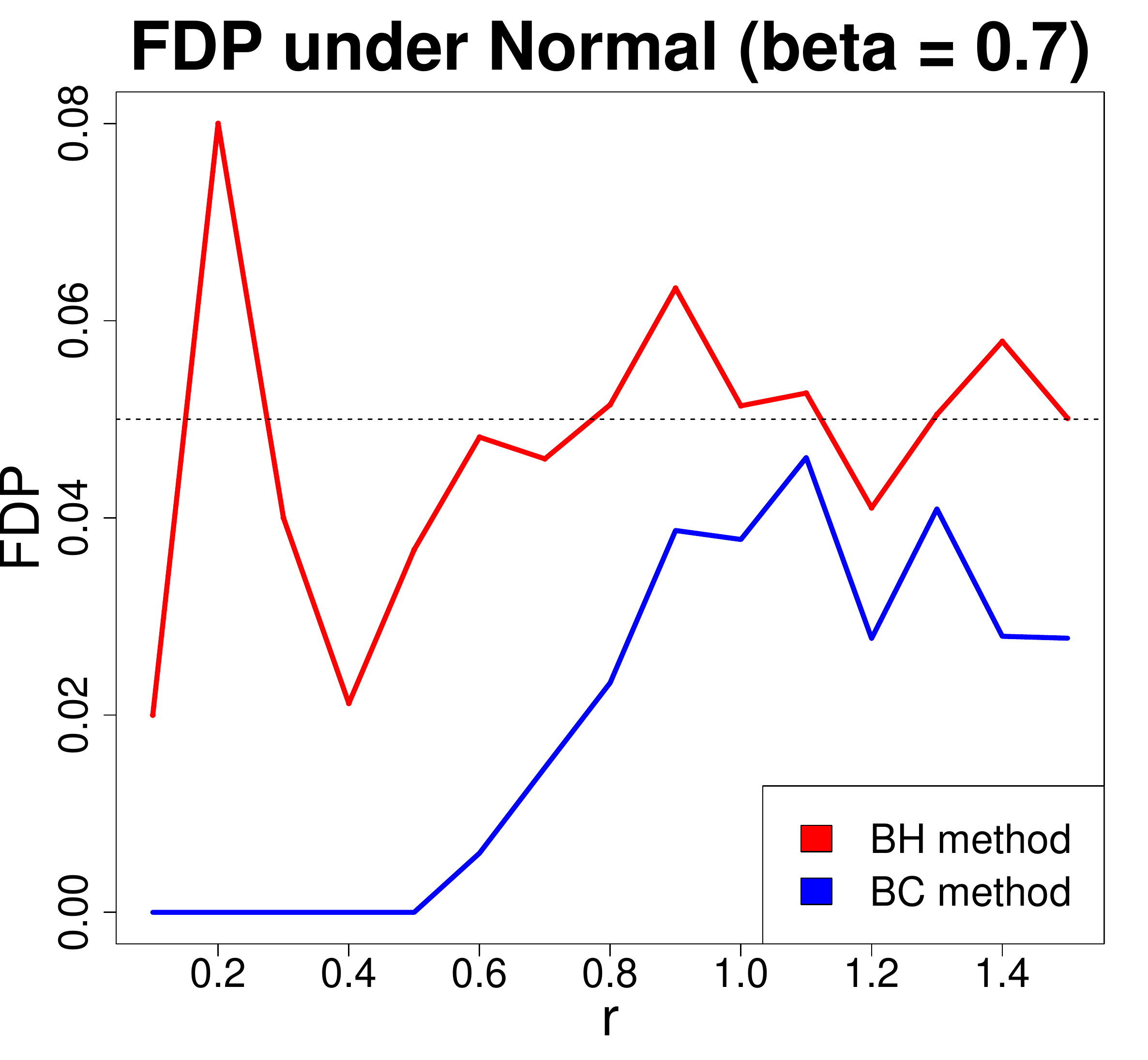}
	\caption{Simulation results showing the FDP for the BH and BC methods under the normal model in three distinct sparsity regimes. The black horizontal line delineates the desired FDR control level ($q = 0.05$).}
	\label{fig:fdp_normal}	
\end{figure}

\begin{figure}[h!]\centering
	\includegraphics[width = 5cm, height= 5cm]{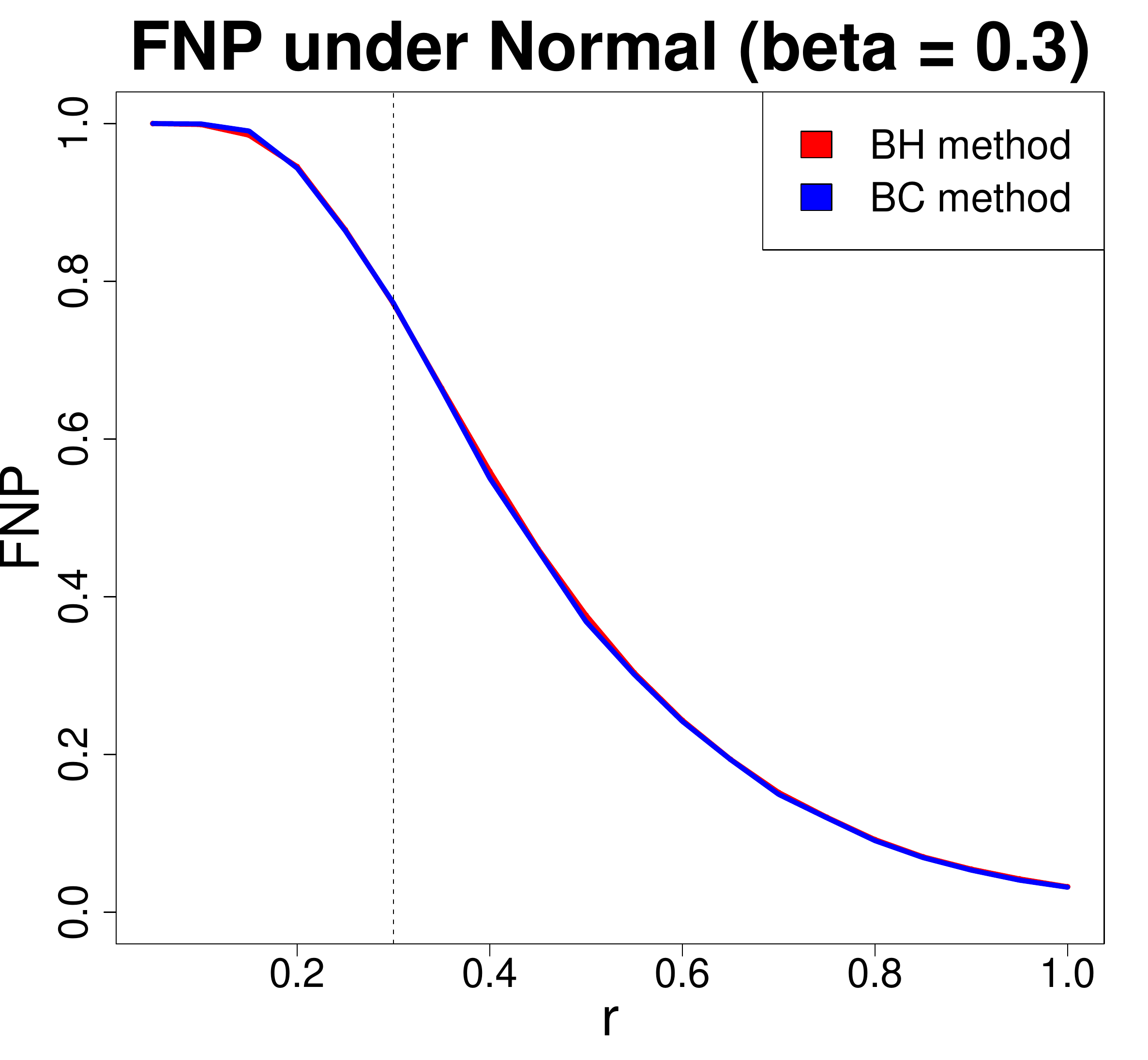}
	\includegraphics[width = 5cm, height= 5cm]{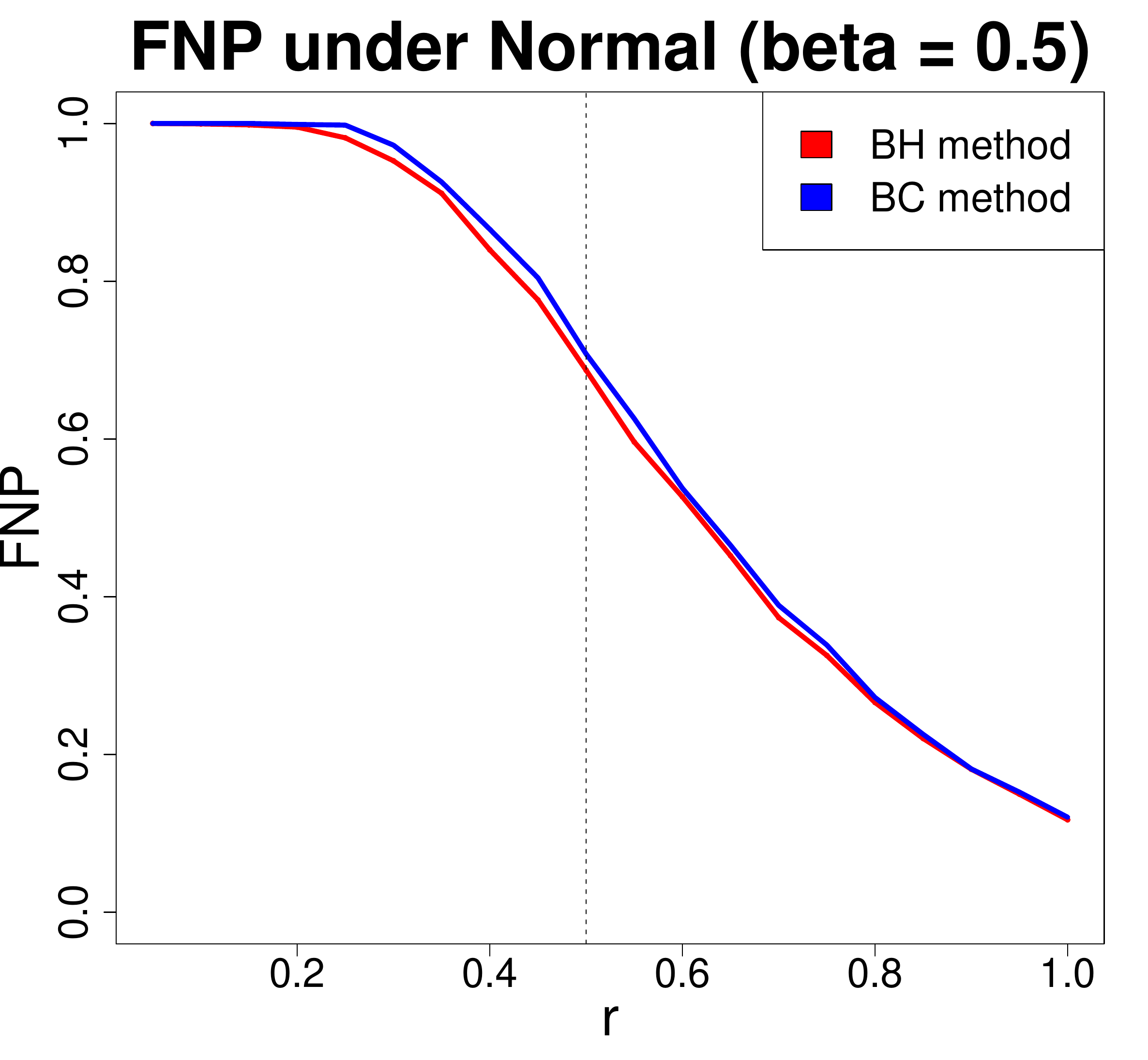}
	\includegraphics[width = 5cm, height= 5cm]{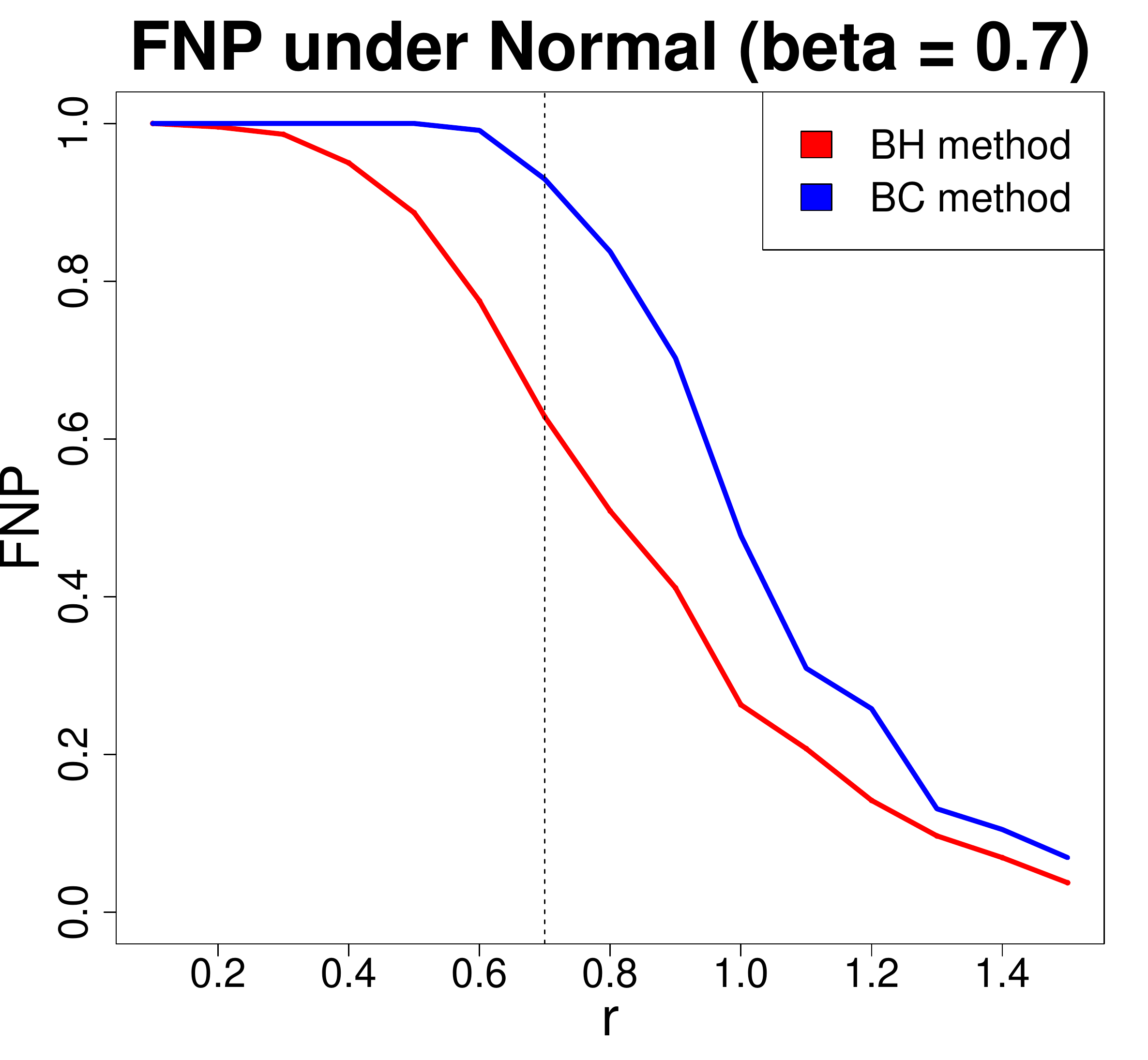}
	\caption{Simulation results showing the false non-discovery proportion (FNP) for the BH and BC methods under the normal model in three distinct sparsity regimes. The black vertical line delineates the theoretical threshold  ($r=\beta$).}
	\label{fig:fnp_normal}
\end{figure}

In this model, $\Psi$ is the standard normal distribution.  
The simulation results are reported in \figref{fdp_normal} and \figref{fnp_normal}.
In \figref{fdp_normal} we report the FDP.  Recall that the methods are set to control the FDR at the desired level ($q = 0.05$).  We see that the BC method becomes more conservative than the BH method as $\beta$ increases.  
In \figref{fnp_normal} we report the FNP.  
We see that the BC method performs comparably to the (oracle) BH method  at $\beta = 0.3$ and $\beta = 0.5$, but is clearly less powerful in the sparsest regime $\beta = 0.7$.  This is in line with the earlier observation that the BC method becomes more conservative with increasing values of $\beta$.  It can also be explained by the fact, at $\beta = 0.7$, the number of false nulls ($m = 31$ out of $n = 10^5$) is too small to reveal the asymptotic power of the BC method.  
Finally, we remark that the transition from high FNP to low FNP happens in the vicinity of the theoretical threshold ($r = \beta$). 

\subsubsection{Double-exponential model}

\begin{figure}[h!]\centering
	\includegraphics[width = 5cm, height= 5cm]{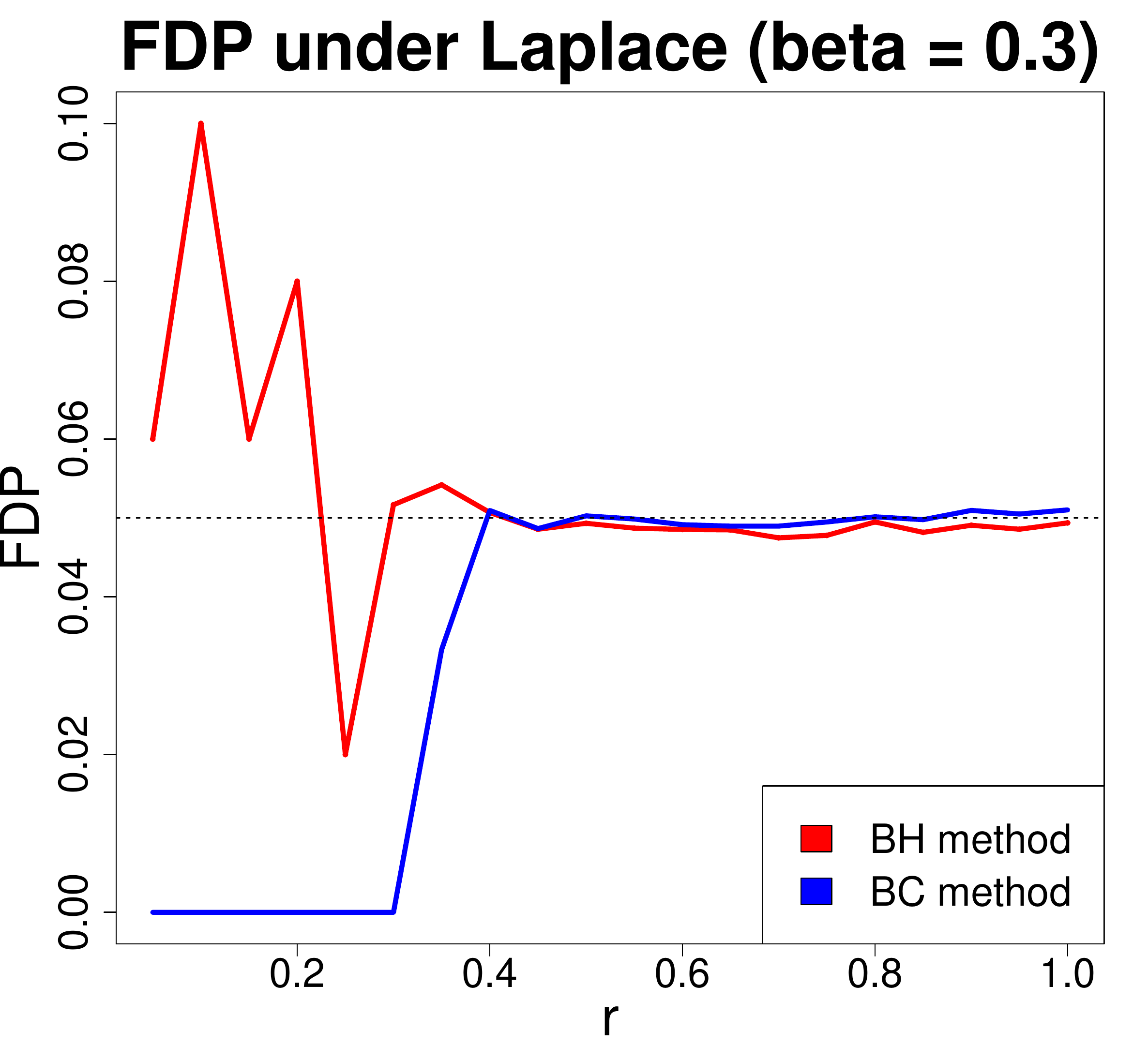}
	\includegraphics[width = 5cm, height= 5cm]{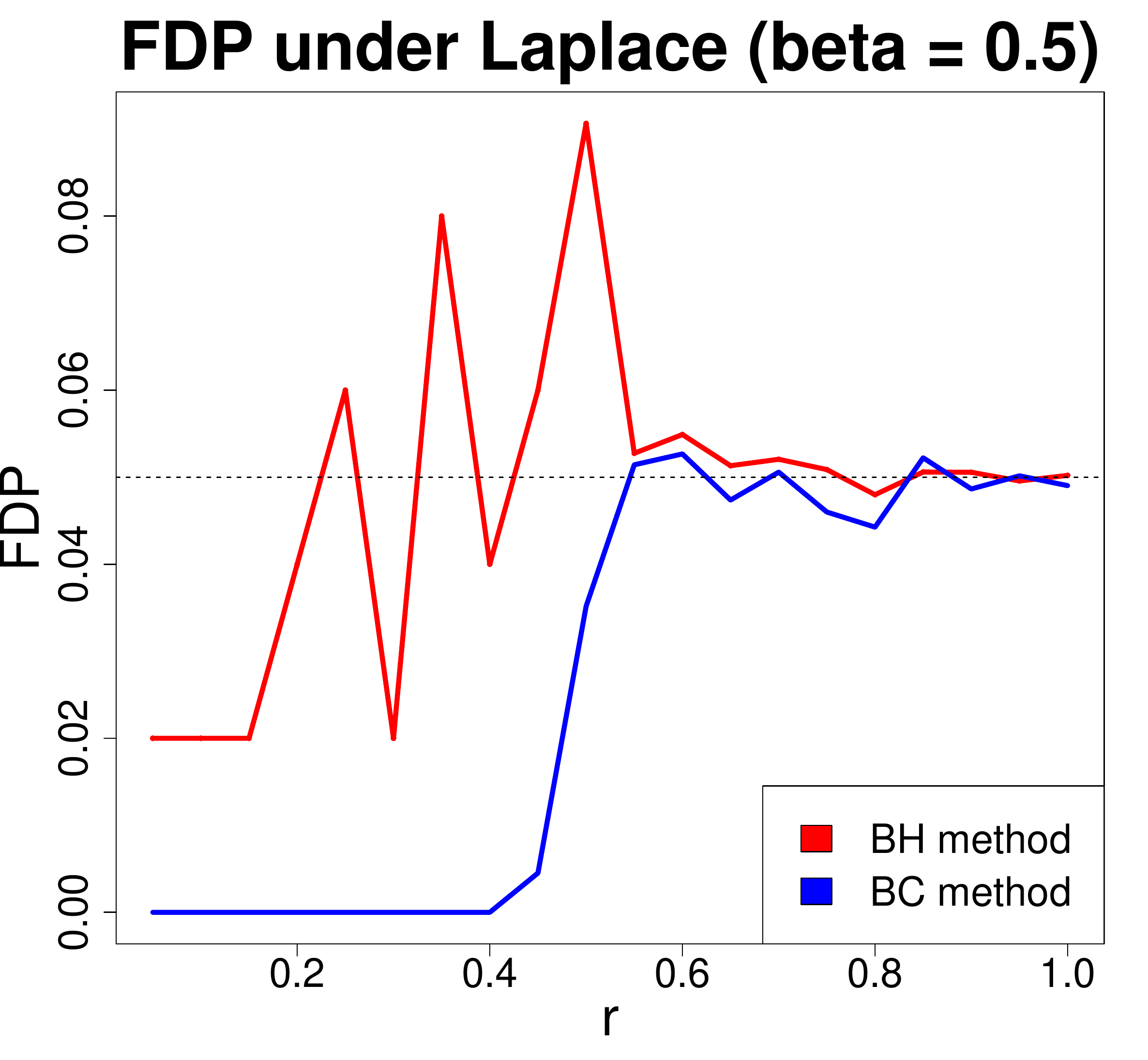}
	\includegraphics[width = 5cm, height= 5cm]{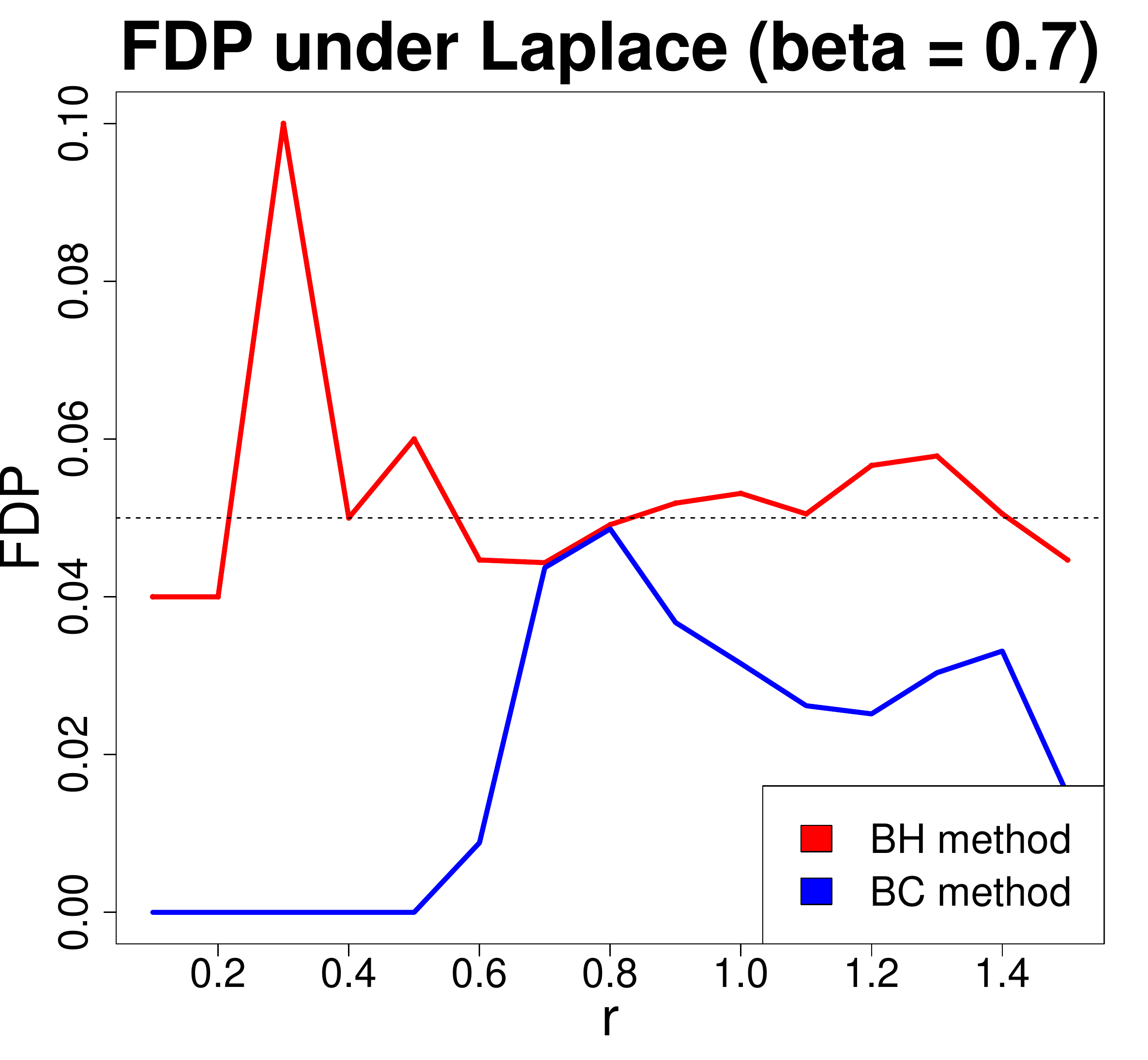}
	\caption{Simulation results of the false discovery proportion (FDP) for double-exponential model in three distinct sparsity regimes. The black horizontal line delineates the desired FDR control level ($q = 0.05$).}
	\label{fig:fdp_laplace}	
\end{figure}

\begin{figure}[h!]\centering
	\includegraphics[width = 5cm, height= 5cm]{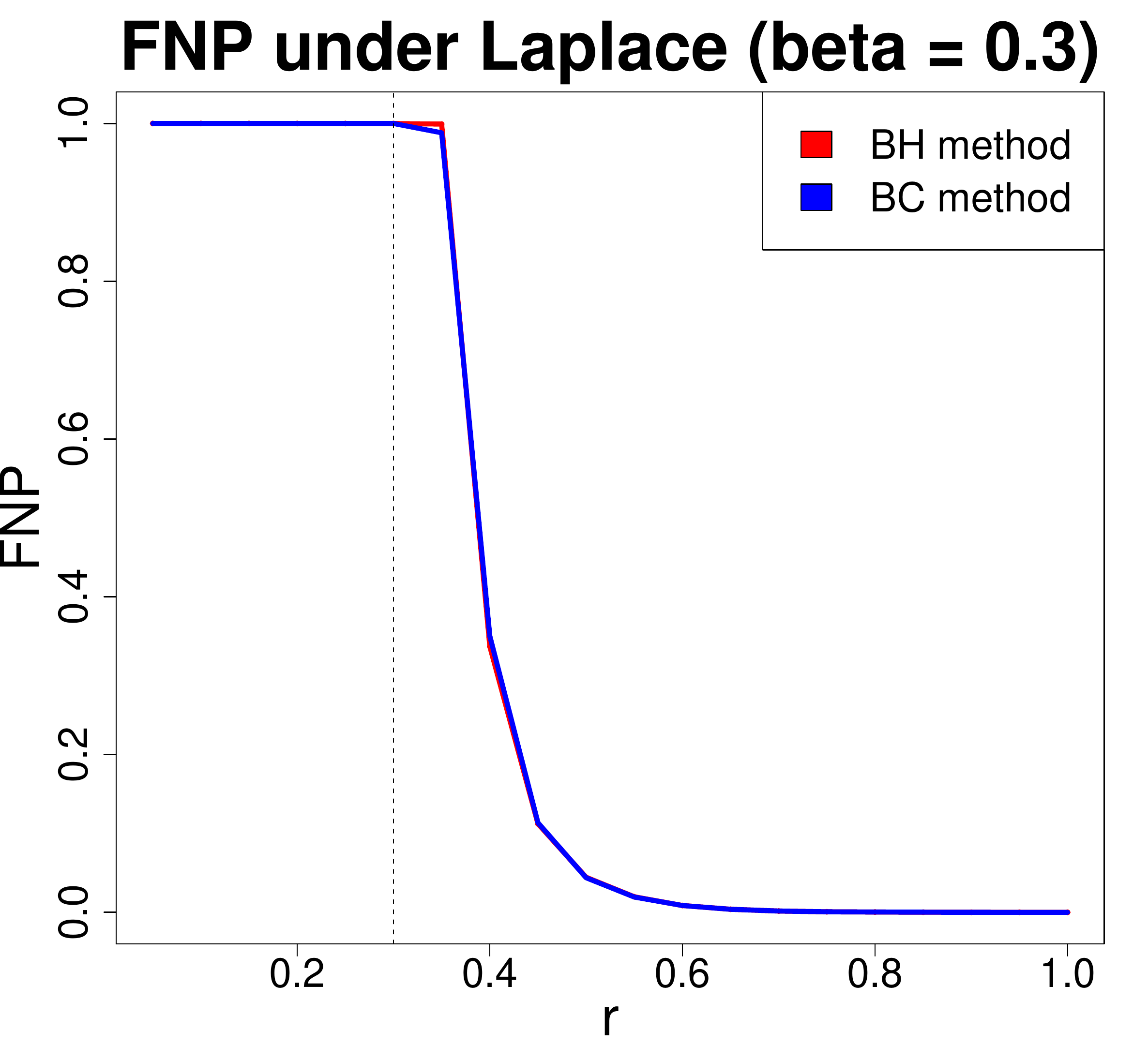}
	\includegraphics[width = 5cm, height= 5cm]{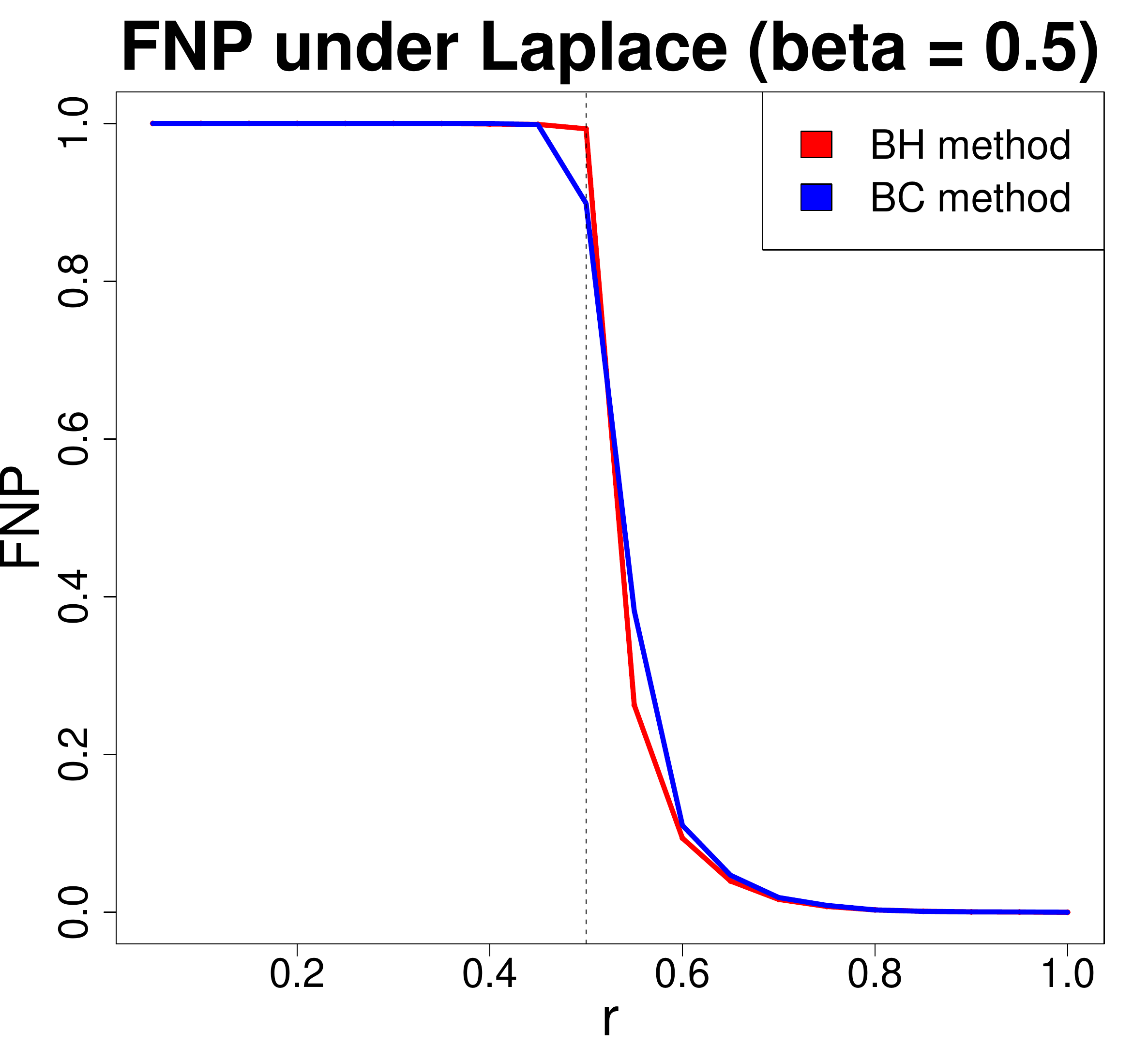}
	\includegraphics[width = 5cm, height= 5cm]{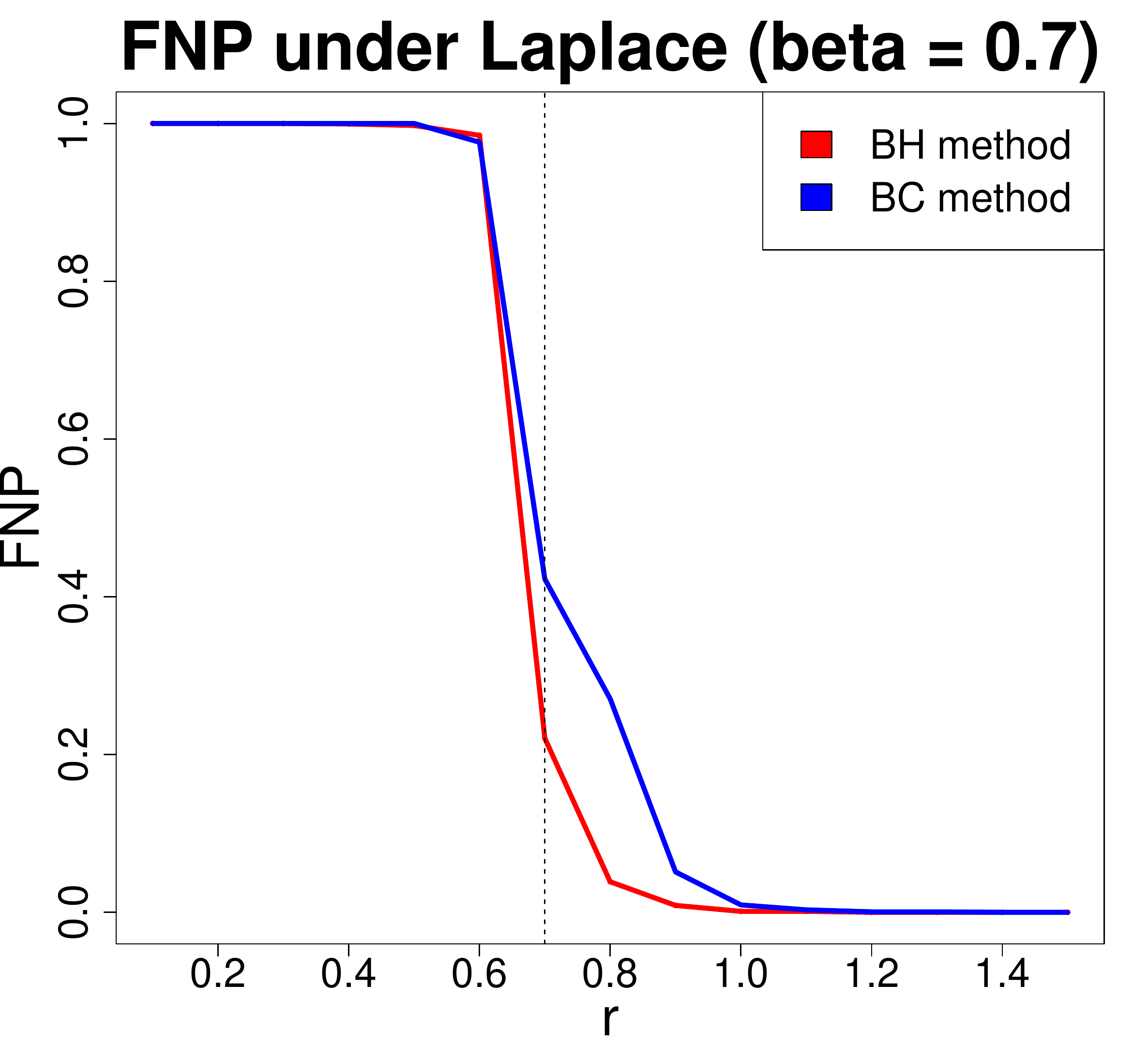}
	\caption{Simulation results of the false non-discovery proportion (FNP) for double-exponential model in three distinct sparsity regimes. The black vertical line delineates the detection boundary ($r=\beta$).}
	\label{fig:fnp_laplace}
\end{figure}

In this model, $\Psi$ is double-exponential distribution with variance of 1.  The simulation results are reported in \figref{fdp_laplace} (FDP) and \figref{fnp_laplace} (FNP).  Here we observe that the BC method is rather conservative regardless of $\beta$.  The two methods are again comparable in terms of FNP, in fact a bit more so than in the normal setting.  The transition from FNP near 1 to FNP near 0 happens, again, in the vicinity of the theoretical threshold, but is much sharper here.

\subsection{Varying sample size}

In this second set of experiments, we examine the effect of various sample sizes on the risk of BH and BC procedures under standard normal model and double-exponential model (with variance~1).  We simultaneously explore the effect of letting the desired FDR control level $q$ tend to 0, in accordance with \eqref{q}.  Specifically, we set it as $q = q_n = 1/\log n$.  We choose $n$ on a log scale, specifically, $n \in \{10^2, 10^3, 10^4, 10^5, 10^6\}$.  Each time, we fix a value of $(\beta, r)$ such that $r > \beta$. 

\begin{figure}[h!]\centering
		\includegraphics[width = 5cm, height= 5cm]{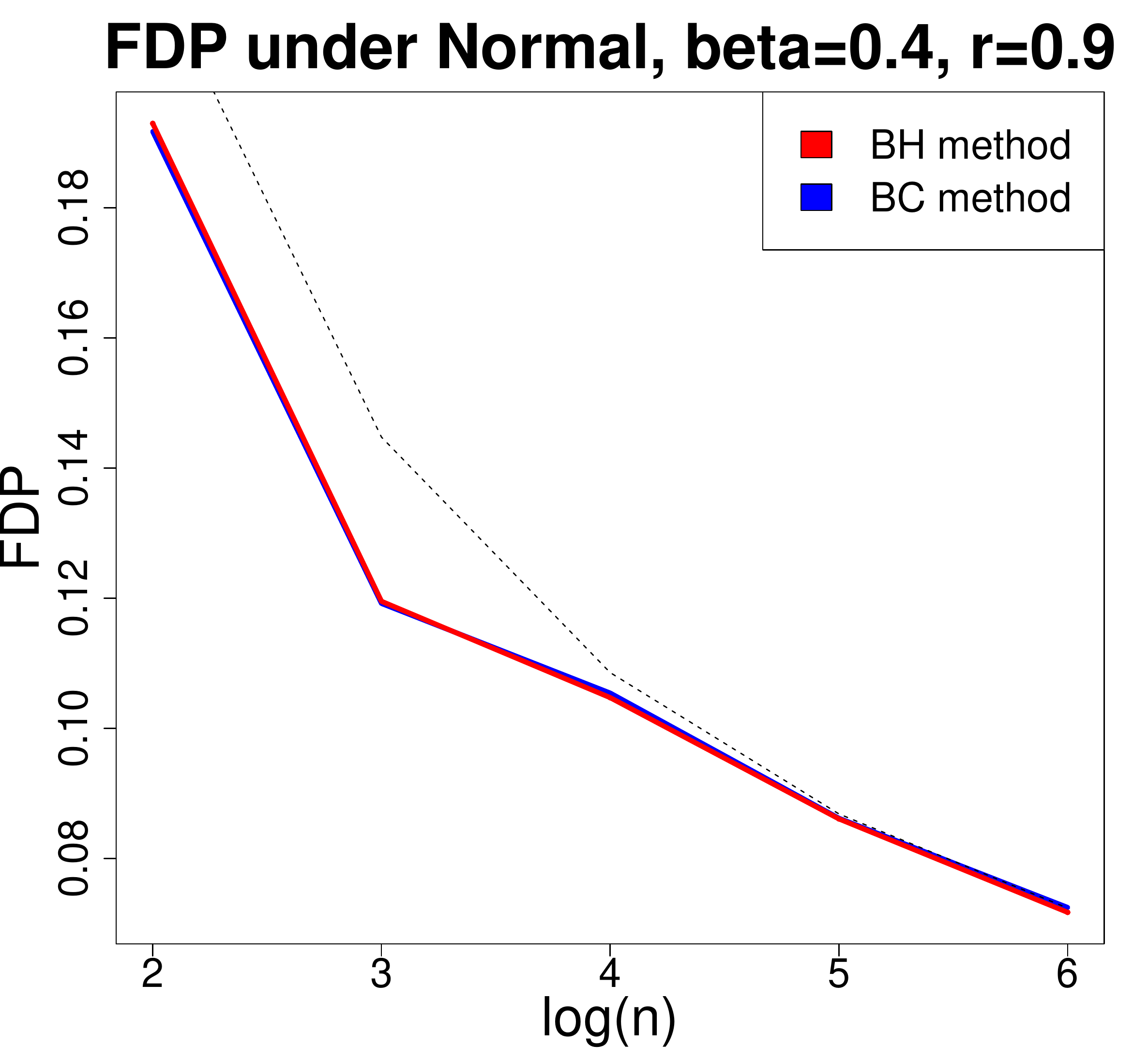}
		\includegraphics[width = 5cm, height=5cm]{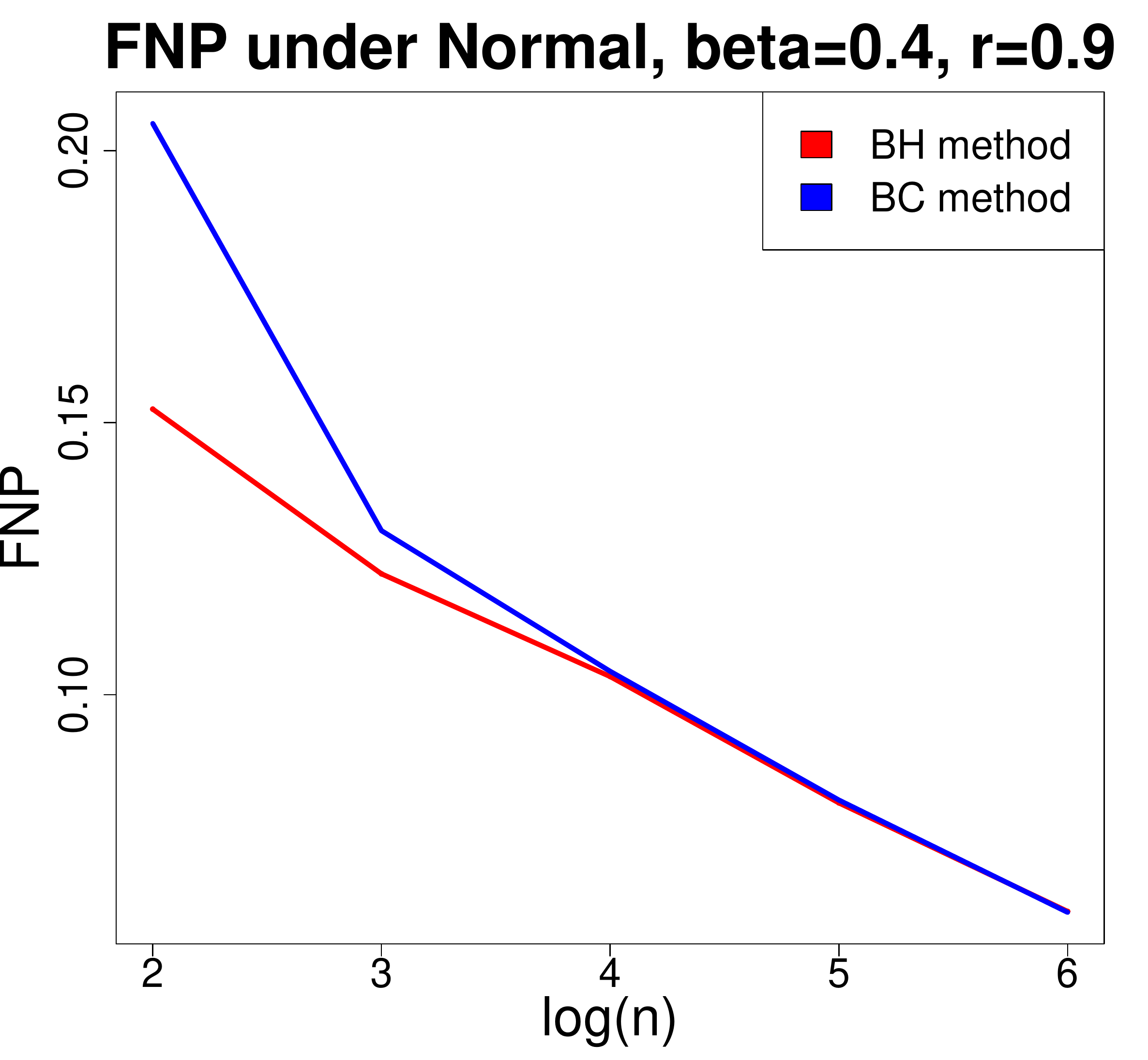}
		\caption{FDP and FNP in the normal model under $(\beta, r) = (0.4, 0.9)$ with varying sample size.}
		\label{fig:vary_dense_normal}
\end{figure}

\begin{figure}[h!]\centering
		\includegraphics[width = 5cm, height= 5cm]{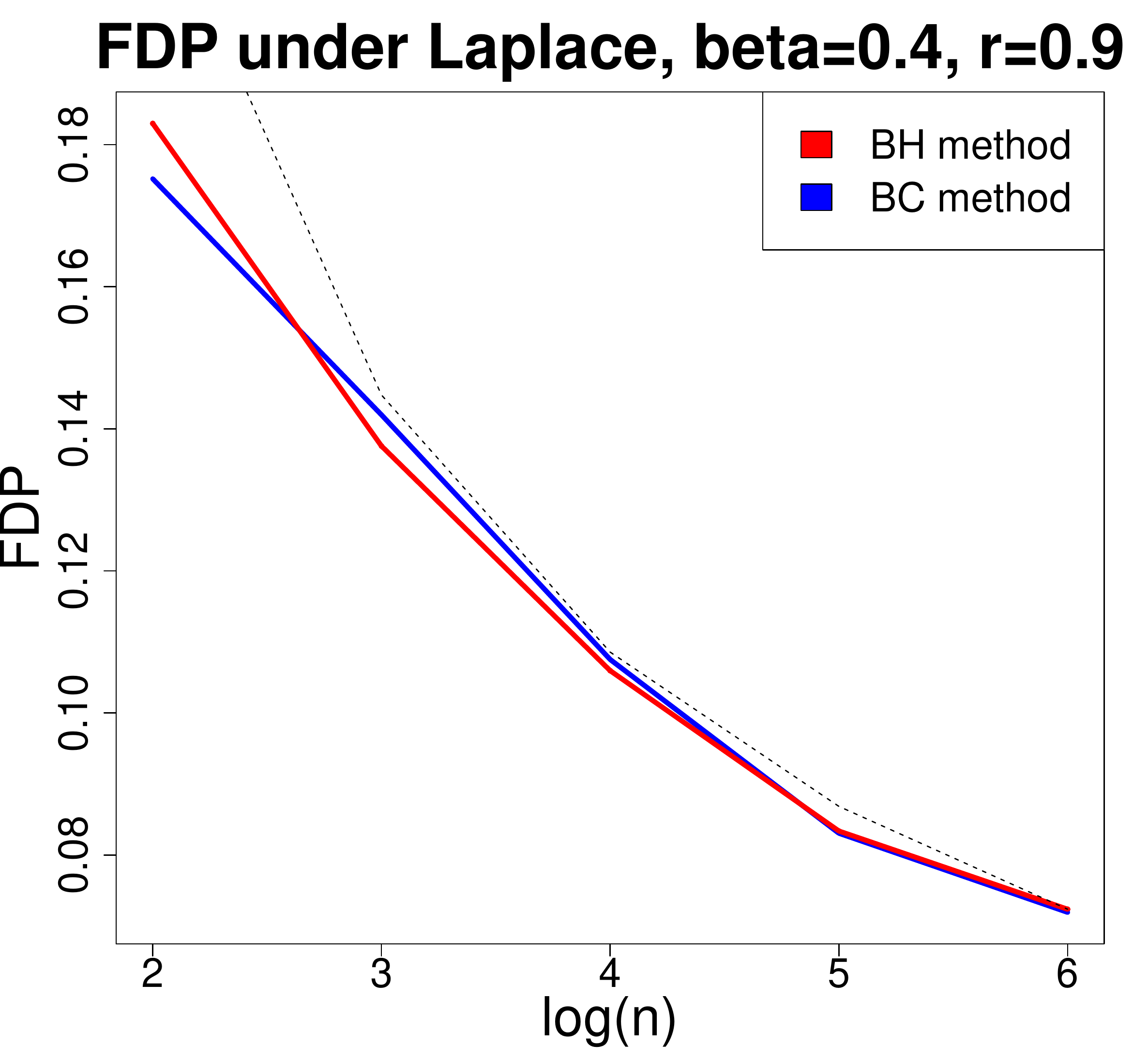}
		\includegraphics[width = 5cm, height= 5cm]{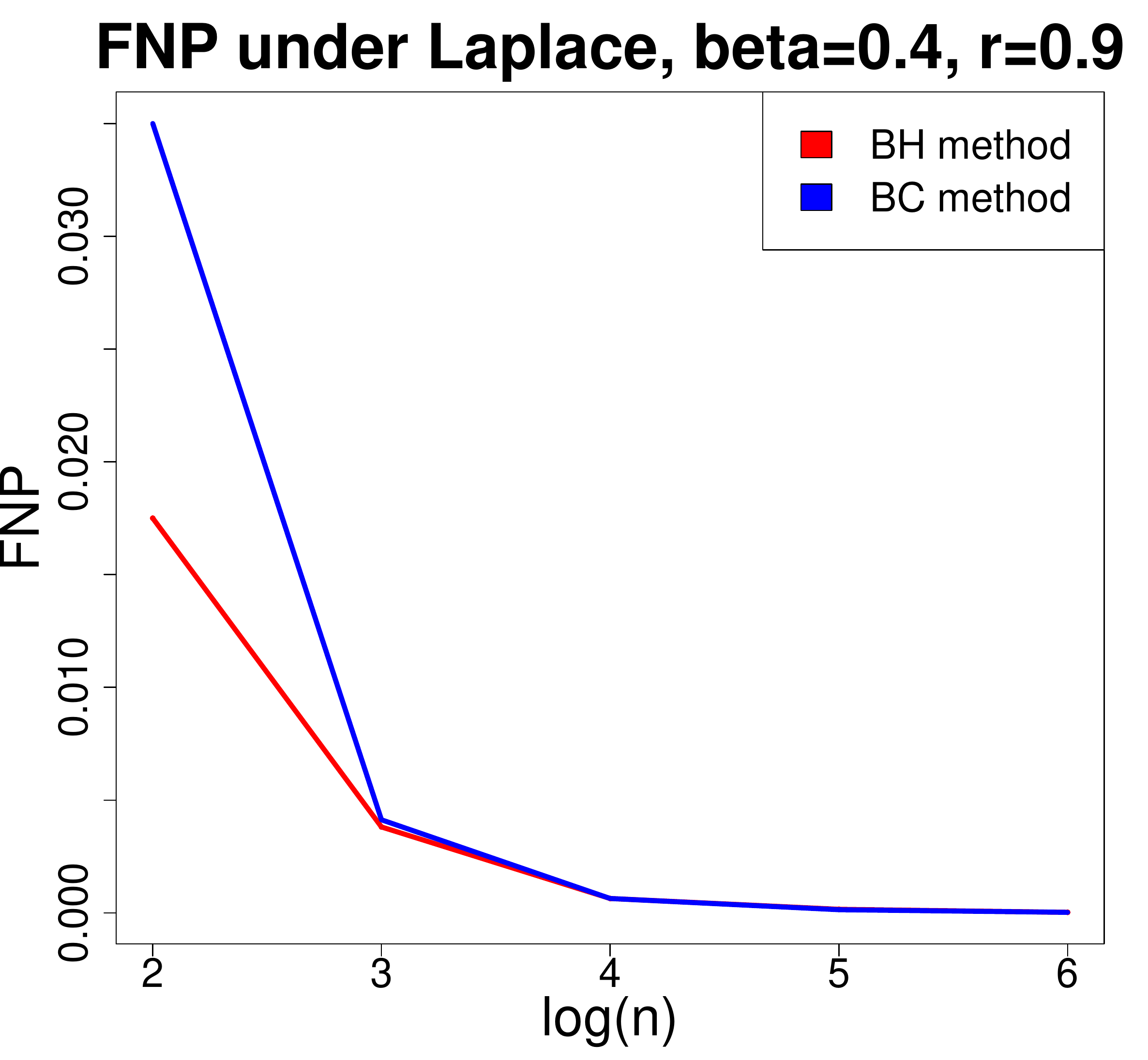}
		\caption{FDP and FNP in the double-exponential model under $(\beta, r) = (0.4, 0.9)$ with varying sample size.}
		\label{fig:vary_dense_laplace}
\end{figure}

In the first setting, we set $(\beta, r) = (0.4, 0.9)$.
The simulation results are reported in \figref{vary_dense_normal} and \figref{vary_dense_laplace}.
We see that, in both models, the risks of the two procedures decrease to zero rapidly as the sample size gets larger.  The BH method clearly dominates (in terms of FNP) up until $n = 10^3$, and after that the two methods behave similarly.

\begin{figure}[h!]\centering
		\includegraphics[width = 5cm, height= 5cm]{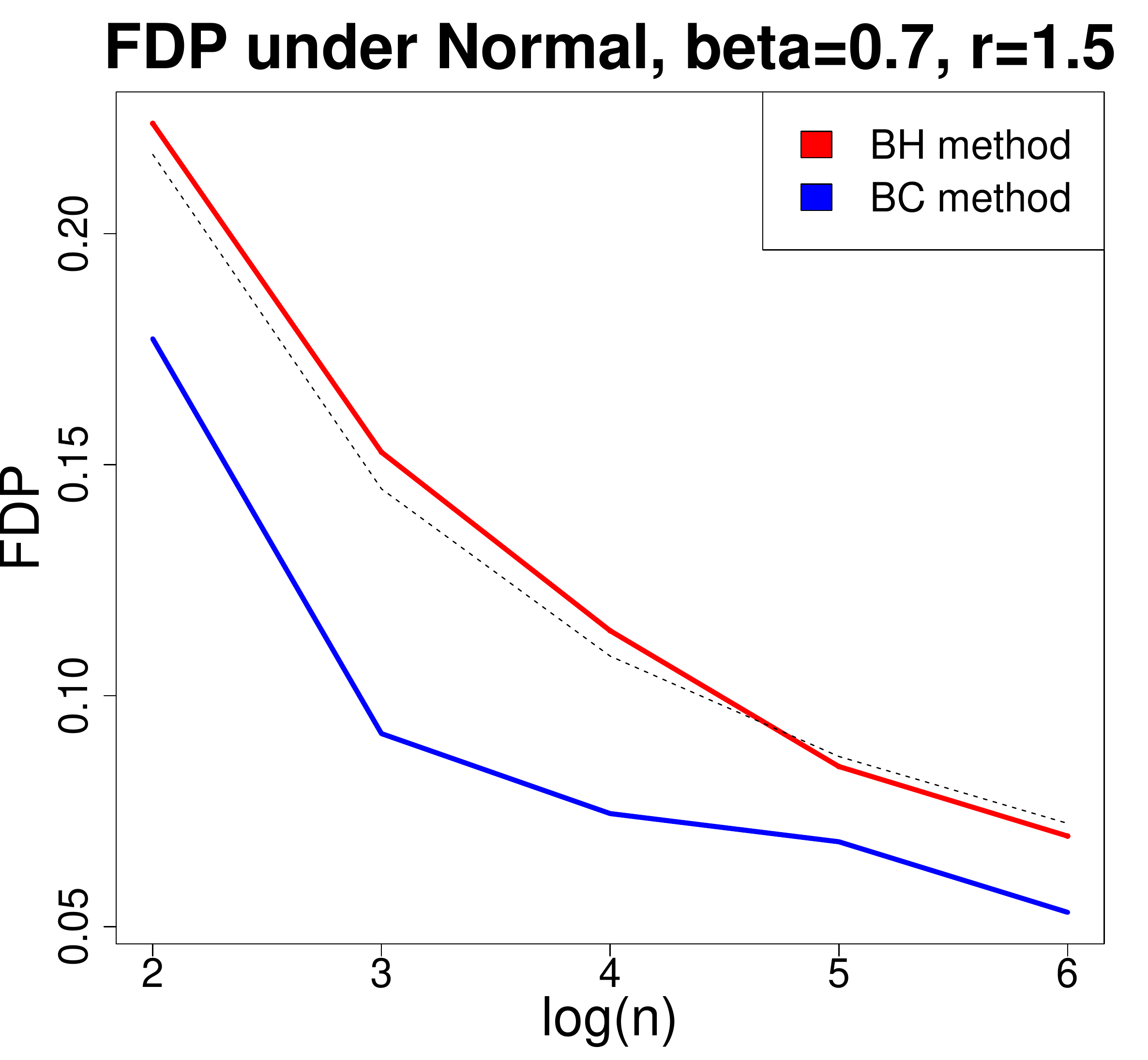}
		\includegraphics[width = 5cm, height= 5cm]{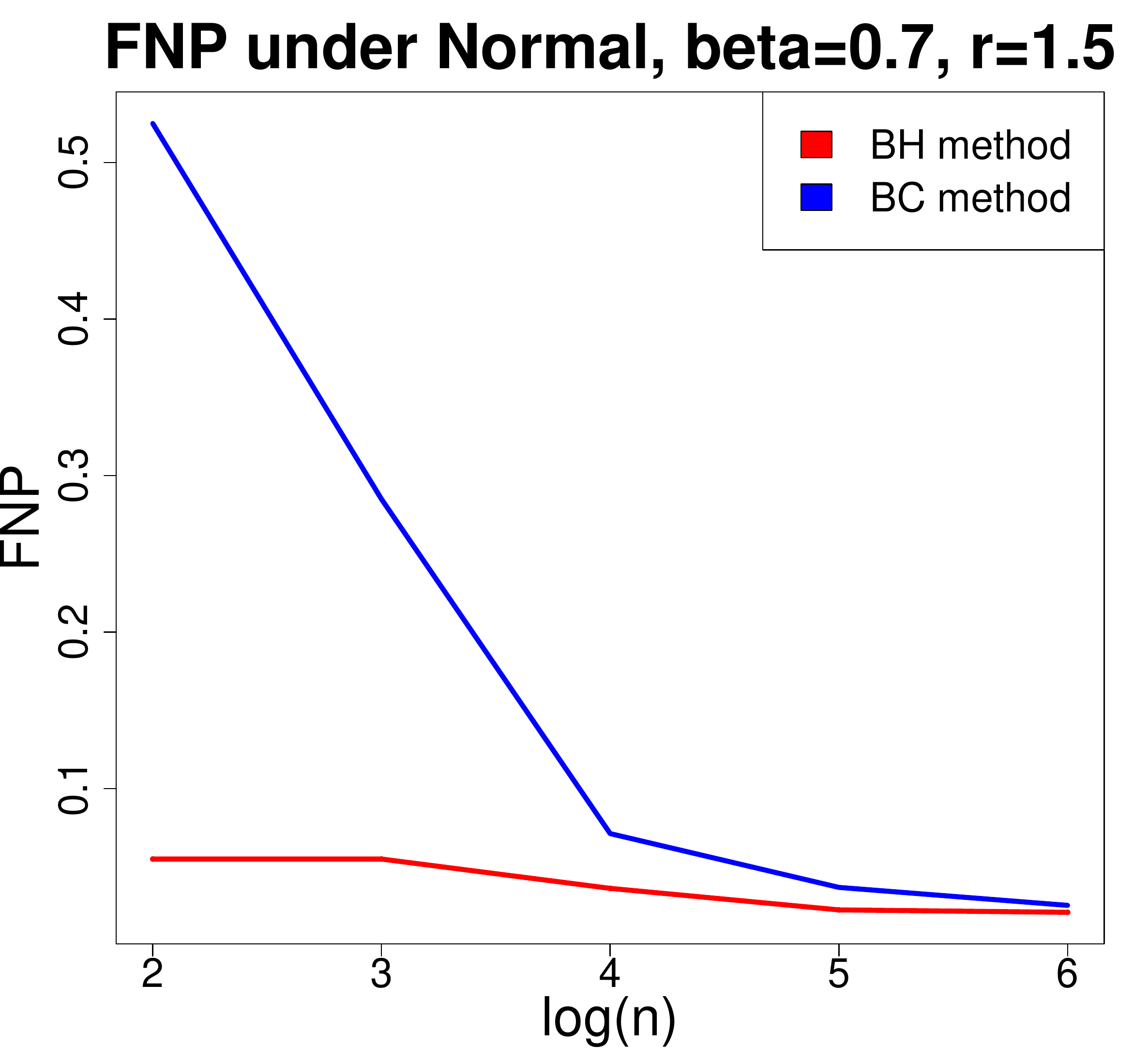}
		\caption{FDP and FNP in the normal model under $(\beta, r) = (0.7, 1.5)$ with varying sample size.}
		\label{fig:vary_sparse_normal}
\end{figure}

\begin{figure}[h!]\centering
		\includegraphics[width = 5cm, height= 5cm]{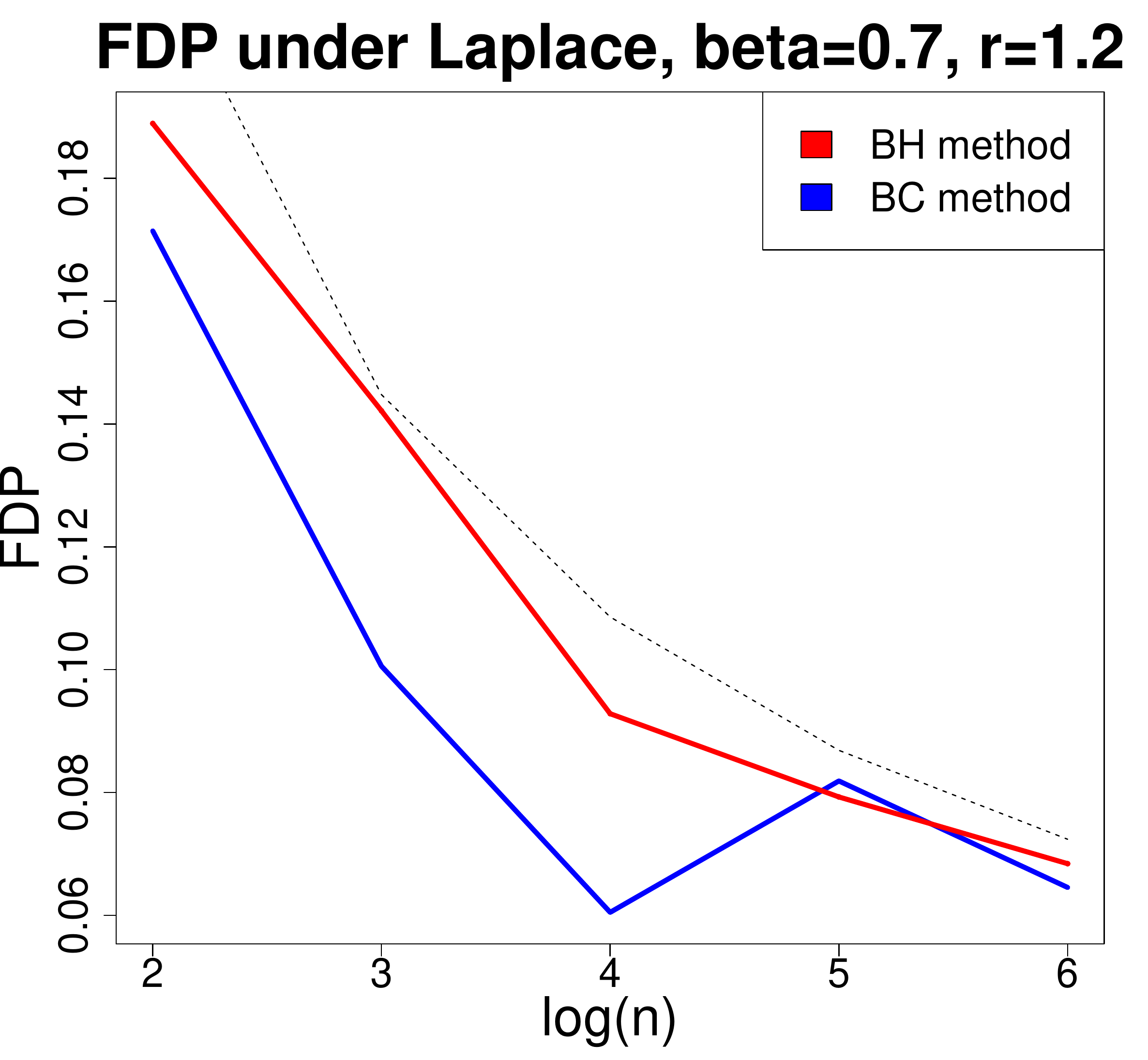}
		\includegraphics[width = 5cm, height= 5cm]{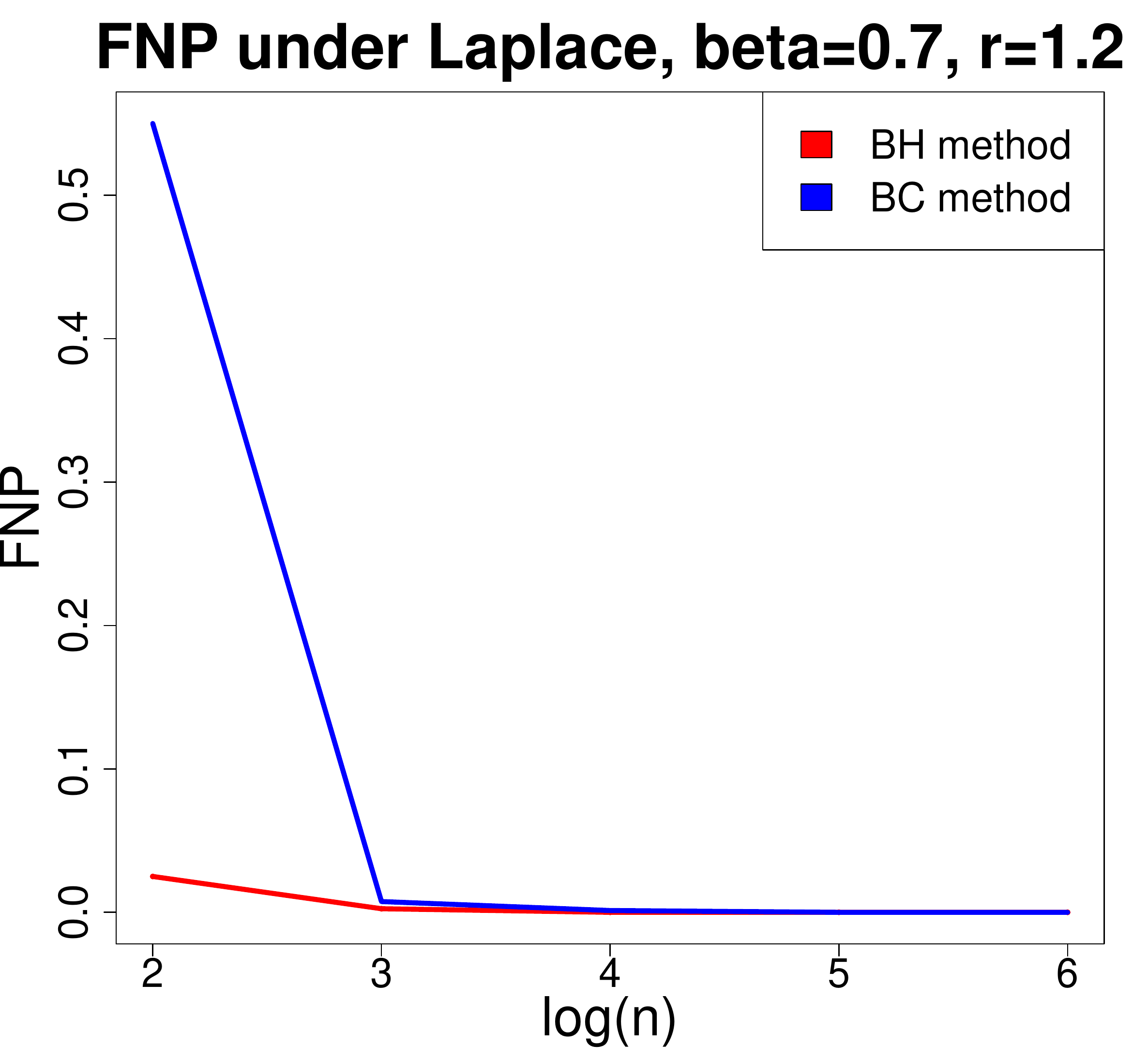}
		\caption{FDP and FNP in the double-exponential model under $(\beta, r) = (0.7, 1.2)$ with varying sample size.}
		\label{fig:vary_sparse_laplace}
\end{figure}

In the second setting, we set $(\beta, r) = (0.7, 1.5)$ for normal model and $(\beta, r) = (0.7, 1.2)$ for double-exponential model. 
The simulation results are reported in \figref{vary_sparse_normal} and \figref{vary_sparse_laplace}.
In this sparser regime, we can see that the BC method is much more conservative than BH method when $n$ is relatively small. But as $n$ gets larger, this is less pronounced.  The BH method clearly dominates (in terms of FNP) up until $n = 10^3$ and past $n = 10^4$ the two methods behave similarly.  The difference is much more dramatic here, in line with our findings in \secref{fixed}.

\section{Proofs}
\label{sec:proofs}

We prove our results in this section.

%

\subsection{Proof of \thmref{lower}} \label{sec:proof-thm-lower}
We first remark that, for any procedure (left implicit)
\beq
\fdp + \fnp \ge 1 - \frac{|\F|}{|\R|} \wedge \frac{|\R|}{|\F|}.
\eeq
Hence, to show that the procedure has risk tending to 1, by dominated converge, it suffices to show that $|\F|/|\R| \to 0$ or $|\R|/|\F| \to 0$.  
This brutish tactic is in fact enough in all cases below except for the last one.

We have the following simple facts.
On the one hand, 
\beq
|\F| = m \sim n^{1-\beta} = o(n),
\eeq 
under \eqref{eps}. 
On the other hand, for any $t \in \bbR$, we have
\beq
|\R_t \cap \F| \sim \Bin(m, \Psi(t -\mu)) \quad \indep \quad
|\R_t \setminus \F| \sim \Bin(n-m, \Psi(t)).
\eeq
In particular, by Chebyshev's inequality,
\begin{align}
|\R_t \cap \F| &= m \Psi(t -\mu) + O_\bbP(m \Psi(t -\mu))^{1/2}, \label{RcapF} \\
|\R_t \setminus \F| &= (n-m) \Psi(t) + O_\bbP((n-m) \Psi(t))^{1/2}. \label{RminusF}
\end{align}
We will also use the fact that 
\beq \label{Rsum}
|\R_t| = |\R_t \cap \F| + |\R_t \setminus \F|, \quad
|\F \setminus \R| = |\F| - |\R_t \cap \F|.
\eeq

By shifting $\Psi$ if needed, we may assume with loss of generality that $\Psi(0) = 1/2$.  

\paragraph{Case: $t \le 0$.}  
In this case there are too many rejections.  Indeed, for any $t \le 0$, 
\beq
|\R_t| \ge |\R_0| = (1/2 + o_\bbP(1)) n \gg m = |\F|,
\eeq
so that $|\F|/|\R_t| \to 0$.

\paragraph{Case: $t > 0$.}  
Such a threshold $t > 0$ is necessarily of the form $t = (\gamma r_n \log n)^{1/\gamma}$ for some $r_n > 0$. 
Extracting subsequences if needed, we may restrict ourselves to the following situations.
\bitem
\item {\bf Subcase: $r_n \to r_\infty \in (\beta, \infty]$.}  
In this case there are not enough rejections.  Indeed, on the one hand, we have \eqref{RcapF}, and using \defref{AGG}, $\log \Psi(t) \sim - t^\gamma/\gamma = - r_n \log n$, we have $n \Psi(t) = n^{1-r_n+o(1)} \ll n^{1-\beta} \sim m$, because $r_n \to r_\infty > \beta$.
On the other hand, we have \eqref{RminusF}, and using the fact that $t - \mu \to \infty$, we have $\Psi(t-\mu) \to 0$, and in particular, $m \Psi(t-\mu)  = o(m)$.
This and \eqref{Rsum} implies that $|\R_t|/|\F| \to 0$.

\item {\bf Subcase: $r_n \to r_\infty \in [0,\beta)$.}  
In this case there are too many rejections.  Indeed, $|\R_t| \ge |\R_t \setminus \F| = (1 + o_\bbP(1)) n \Psi(t)$, with $n\Psi(t) = n^{1-r_n + o(1)} \gg n^{1-\beta} \sim m$, because $r_n \to r_\infty < \beta$.  This implies that $|\F|/|\R_t| \to 0$.

\item {\bf Subcase: $r_n \to r_\infty = \beta$.}  
In this case, the number of rejections is about right, but the rejections themselves are not accurate.  Indeed, we have 
\beq
|\R_t \cap \F| = m \Psi(t -\mu) + O_\bbP(m \Psi(t -\mu))^{1/2},
\eeq
with $\Psi(t -\mu) = o(1)$ since $t -\mu \to \infty$ due to the fact that $r_n \to\beta > r$.  Hence, $|\R_t \cap \F|/|\F| \to 0$, which in turn implies that $|\F \setminus \R_t|/|\F| \to 1$.
\eitem

\subsection{Proof of \thmref{BH}}
\label{sec:proof-BH}

Let $\Psi$ denote the null survival function, assumed to be AGG with parameter $\gamma \ge 1$.  
Let $\hat G$ denote the empirical survival function
\beq
\hat G(t) = \frac1n \sum_{i \in [n]} \IND{X_i \ge t}.
\eeq
Let $Y_i = X_i - \mu_i$ and note that $(Y_i : i \in [n])$ are IID with distribution $\Psi$.  Define the empirical survival functions 
\beq
\hat W_{\rm true}(y) = \frac1{n-m} \sum_{i \notin \cF} \IND{Y_i \ge y},  \quad
\hat W_{\rm false}(y) = \frac1m \sum_{i \in \cF} \IND{Y_i \ge y},
\eeq
so that
\beq
\hat G(t) = (1-\eps) \hat W_{\rm true}(t) + \eps \hat W_{\rm false}(t - \mu).
\eeq
where $\eps := m/n \sim n^{-\beta}$ under \eqref{eps}.  

We need the following result to control the deviations of the empirical distributions. 

\begin{lem}[\cite{eicker1979asymptotic}] \label{lem:emp}
Let $Z_1, \dots, Z_k$ be IID with continuous survival function $Q$.  Let $\hat Q_k$ denote their empirical survival function and define $\zeta_k = \sqrt{2 \log\log(k)/k}$ for $k \ge 3$.  Then
\beq
\frac1{\zeta_k} \max_z \frac{\hat Q_k(z) - Q(z)}{\sqrt{Q(z) (1-Q(z))}} \to 1, \text{ in probability as } k \to \infty.
\eeq
In particular, 
\beq
\hat Q_k(z) = Q(z) + O_\bbP(\zeta_k) \sqrt{Q(z) (1-Q(z))}, \quad \text{uniformly in $z$}.
\eeq
\end{lem}

Applying \lemref{emp}, we get
\begin{align}
\hat G(t) 
&= (1-\eps) \big[\Psi(t) + O_\bbP(\zeta_n) \sqrt{\Psi(t) (1-\Psi(t))}\big] \\
&\quad + \eps \big[\Psi(t-\mu) + O_\bbP(\zeta_m) \sqrt{\Psi(t-\mu) (1-\Psi(t-\mu))}\big].
\end{align}
From this we get
\beq \label{Ghat}
\hat G(t) = G(t) + \hat R(t),
\eeq 
where
\beq\label{G}
G(t) := \E[\hat G(t)] = (1-\eps) \Psi(t) + \eps \Psi(t-\mu),
\eeq
and
\beq\label{Rhat}
\hat R(t) = O_\bbP\big(\zeta_n \sqrt{\Psi(t) (1-\Psi(t))} + \zeta_m \eps \sqrt{\Psi(t-\mu) (1-\Psi(t-\mu))}\big), \quad \text{uniformly in $t \in \bbR$}.
\eeq

Let $\iota = \iota_{\rm BH}$ be defined as in \eqref{def_BH}.
We have $\hat G(X_{(i)}) = i/n$, so that 
\beq
P_{(i)} \leq i q/n \text{ if and only if } X_{(i)} \ge \Psi^{-1}(q \hat G(X_{(i)})).
\eeq
We therefore have $X_{(i)} \ge \Psi^{-1}(q \hat G(X_{(i)}))$ for $i \le \iota$ and $X_{(i)} < \Psi^{-1}(q \hat G(X_{(i)}))$ for $i > \iota$.
Based on that, and the fact that $\hat G$ is constant between two consecutive $X_i$'s, we have that there is $\tau \in [X_{(\iota)}, X_{(\iota+1)})$ such that  
\beq \label{rho}
\tau = \min \big\{t: t \ge \Psi^{-1}(q\hat G(t))\big\} = \min \big\{t: t = \Psi^{-1}(q\hat G(t))\big\}.
\eeq
Note that the BH procedure coincides with $\R_\tau$, the threshold method with threshold $\tau$.  In particular,
\beq\label{Fhat}
\fnp(\R_\tau) = 1 - \hat F(\tau), \quad 
\hat F(t) := \frac1m \sum_{i \in \F} \IND{X_i \ge t},
\eeq
so that it suffices to show that $\hat F(\tau) \to 1$ in probability.
As above, by \lemref{emp},
\beq
\hat F(t) = \hat W_{\rm false}(t - \mu) = \Psi(t - \mu) + O_\bbP(\zeta_m) \sqrt{\Psi(t-\mu) (1-\Psi(t-\mu))},
\eeq
and in particular $\hat F(\tau) = \Psi(\tau-\mu) + o_\bbP(1)$, so it suffices to show that $\tau-\mu \to -\infty$ in probability.

Since $r>\beta$, we can take a real number $r_* \in (\beta, r)$.
Define $t_* = (\gamma r_* \log n)^{1/\gamma}$.  Since $t_* - \mu \to -\infty$, it suffices to show that $\tau \le t_*$ with probability tending to 1. 
We have 
\beq
G(t_*) = (1-\eps) \Psi(t_*) + \eps \Psi(t_* - \mu).
\eeq 
The first term is $\sim \Psi(t_*)$, with 
\beq\label{Psi-tstar}
\Psi(t_*) = n^{-r_* + o(1)},
\eeq 
by \defref{AGG}, which says that $\log \Psi(t) \sim - t^\gamma/\gamma$ as $t \to \infty$.
The second term is $\sim n^{-\beta}$ by \eqref{eps} and the fact that $\Psi(t_* - \mu) \to 1$ since, again, $t_* - \mu \to -\infty$.
Together, we obtain $G(t_*) \sim n^{-\beta}$, using also the fact that $r_* > \beta$.
In addition, by \eqref{Rhat} we have 
\beq
\hat R(t_*) 
= O_\bbP\big(\zeta_n \sqrt{\Psi(t_*)}) + o_\bbP(\eps) 
= o_\bbP(n^{-\beta}),
\eeq
since $\zeta_n \sqrt{\Psi(t_*)} = n^{-\frac12 (r_*+1) +o(1)}$, again by \eqref{Psi-tstar}, and $\beta < r_* < 1$.
Hence, applying \eqref{Ghat}, we obtain
\begin{align}
\hat G(t_*) 
= G(t_*) + \hat R(t_*) 
\sim_\bbP G(t_*)  
\sim n^{-\beta}. 
\end{align}
Together with \eqref{Psi-tstar}, and using by \eqref{q}, we have 
\beq
\hat G(t_*)/\Psi(t_*) = n^{(r_*-\beta) + o_\bbP(1)} \gg 1/q.
\eeq
This, together with \eqref{rho}, implies that $\tau \le t_*$ with probability tending to 1.

\subsection{Proof of \thmref{BC}}
The proof borrows a number of arguments from \secref{proof-BH}.
We use the same notation and assume as before that the $X_i$'s are distinct.
We assume, in addition, that $\Psi$ is symmetric about~0.  

Define the threshold
\beq
\tau = \inf\big\{t : \fdphat(t) \leq q\big\}.
\eeq
The difference with $\tau_{\rm BC}$ in \eqref{BC} is that the range is not limited to $|\bX|$.  It can be seen that $\tau = |X|_{(\iota_{\rm BC}+1)}$ unless $\iota_{\rm BC} = n$ (the BC procedure rejects all the nulls), in which case $\tau = \tau_{\rm BC}$.
This, in particular, implies
\beq
\fnp(\R_\tau) \le \fnp(\R_{\tau_{\rm BC}}) \le \fnp(\R_\tau) + \tfrac1m.
\eeq
Since in our model $m \to \infty$, it suffices to show that $\fnp(\R_\tau) \to 0$ in probability.  As before, \eqref{Fhat} holds true, so it suffices to show that $\hat F(\tau) \to 1$ in probability.  For that, we saw earlier that it suffices to show that $\tau \le t_*$ with probability tending to 1.

We have
\beq
\fdphat(t_*) = \frac{1 + n(1 - \hat G(-t_*))}{1 \vee n \hat G(t_*)}.
\eeq
We already saw that $\hat G(t_*) \sim n^{-\beta}$, so the denominator above is $\sim n^{1-\beta}$ as $n \to \infty$.
For the numerator, by \eqref{Ghat}, we have
\beq
1 - \hat G(-t_*) = 1 - G(-t_*) - \hat R(-t_*).
\eeq
By \eqref{G}, 
\begin{align}
1 - G(-t_*) 
&= (1-\eps)(1 - \Psi(-t_*)) + \eps (1-\Psi(-t_*-\mu)) \\
&= (1-\eps)\Psi(t_*) + \eps \Psi(t_*+\mu) \quad \text{[by symmetry of $\Psi$]} \\
&\sim \Psi(t_*) = n^{-r_* +o(1)}.  \quad \text{[by \eqref{Psi-tstar}]}
\end{align}
By \eqref{Rhat}, 
\begin{align}
\hat R(-t_*) 
&= O_\bbP\big(\zeta_n \sqrt{1-\Psi(-t_*)} + \zeta_m \eps \sqrt{1-\Psi(-t_*-\mu)}\big) \\
&= O_\bbP(\zeta_n \sqrt{\Psi(t_*)} + \zeta_m \eps \sqrt{\Psi(t_*+\mu)}) \quad \text{[by symmetry of $\Psi$)]} \\
&= O_\bbP(n^{-\frac12(r_*+1)+o(1)} + o(n^{-\frac12(r_*+\beta+1)+o(1)})) \quad \text{[by \eqref{Psi-tstar}]} \\
&= O_\bbP(n^{-\frac12(r_*+1)+o(1)}).
\end{align}
Combined with the fact that $r_* < 1$, we get $1 - \hat G(-t_*) \sim n^{-r_*+o(1)}$, and therefore
\beq
\fdphat(t_*) = \frac{n^{1-r_*+o(1)}}{n^{1-\beta}} = n^{\beta -r_* +o(1)} \ll q. \quad \text{[by \eqref{q} and $\beta < r_\star$]}
\eeq
Hence, $\fdphat(t_*) \le q$ with with probability tending to 1, and when this is the case, $\tau \le t_*$, by definition of $\tau$ above.

\bibliographystyle{chicago}
\bibliography{ref}

\end{document}